# HOEFFDING-ANOVA DECOMPOSITIONS FOR SYMMETRIC STATISTICS OF EXCHANGEABLE OBSERVATIONS

By Giovanni Peccati

*Universités de Paris VI*

Consider a (possibly infinite) exchangeable sequence $\mathbf{X} = \{X_n : 1 \leq n < N\}$, where $N \in \mathbb{N} \cup \{\infty\}$, with values in a Borel space $(A, \mathcal{A})$, and note $\mathbf{X}_n = (X_1, \ldots, X_n)$. We say that $\mathbf{X}$ is *Hoeffding decomposable* if, for each $n$, every square integrable, centered and symmetric statistic based on $\mathbf{X}_n$ can be written as an orthogonal sum of $n$ $U$-statistics with degenerated and symmetric kernels of increasing order. The only two examples of Hoeffding decomposable sequences studied in the literature are i.i.d. random variables and extractions without replacement from a finite population. In the first part of the paper we establish a necessary and sufficient condition for an exchangeable sequence to be Hoeffding decomposable, that is, called *weak independence*. We show that not every exchangeable sequence is weakly independent, and, therefore, that not every exchangeable sequence is Hoeffding decomposable. In the second part we apply our results to a class of exchangeable and weakly independent random vectors $\mathbf{X}_n^{(\alpha,c)} = (X_1^{(\alpha,c)}, \ldots, X_n^{(\alpha,c)})$ whose law is characterized by a positive and finite measure $\alpha(\cdot)$ on $A$ and by a real constant $c$. For instance, if $c = 0$, $\mathbf{X}_n^{(\alpha,c)}$ is a vector of i.i.d. random variables with law $\alpha(\cdot)/\alpha(A)$; if $A$ is finite, $\alpha(\cdot)$ is integer valued and $c = -1$, $\mathbf{X}_n^{(\alpha,c)}$ represents the first $n$ extractions without replacement from a finite population; if $c > 0$, $\mathbf{X}_n^{(\alpha,c)}$ consists of the first $n$ instants of a generalized Pólya urn sequence. For every choice of $\alpha(\cdot)$ and $c$, the Hoeffding-ANOVA decomposition of a symmetric and square integrable statistic $T(\mathbf{X}_n^{(\alpha,c)})$ is explicitly computed in terms of linear combinations of well chosen conditional expectations of $T$. Our formulae generalize and unify the classic results of Hoeffding [*Ann. Math. Statist.* **19** (1948) 293–325] for i.i.d. variables, Zhao and Chen [*Acta Math. Appl. Sinica* **6** (1990) 263–272] and Bloznelis and Götze [*Ann. Statist.* **29** (2001) 353–365 and *Ann. Probab.* **30** (2002) 1238–1265] for finite population statistics. Applications are given to construct infinite "weak urn sequences" and to characterize the covariance of symmetric statistics of generalized urn sequences.

**1. Introduction.** For any $N \in \mathbb{N} \cup \{+\infty\}$, consider a collection $\mathbf{X} = \{X_n : 1 \leq n < N\}$ of exchangeable random observations, whose components take values in some Borel space $(A, \mathcal{A})$ and are defined on a suitable probability space $(\Omega, \mathcal{F}, \mathbb{P})$ [the reader is referred to Aldous (1983) for any unexplained notion concerning exchangeability]. For $1 \leq n < N$ and $q > 0$, we write $\mathbf{X}_n$ and $L^q(\mathbf{X}_n)$, respectively, for the vector $(X_1, \ldots, X_n)$ and for the class of real-valued functionals $T(\mathbf{X}_n)$ such that $\mathbb{E}|T|^q < +\infty$. Roughly speaking, we say that the sequence $\mathbf{X}$ is *Hoeffding decomposable* (or Hoeffding-ANOVA decomposable) if, for every $n$, any centered and symmetric $T \in L^2(\mathbf{X}_n)$ can be uniquely represented as an $L^2$-orthogonal sum of $n$ $U$-statistics based on $\mathbf{X}_n$, say $T_1, \ldots, T_n$, such that each $T_i$ has a (completely) degenerate symmetric kernel of order $i$. In particular, if $\mathbf{X}$ is Hoeffding decomposable, for each $n$ the covariance between symmetric statistics based on $\mathbf{X}_n$ can be represented as a sum of covariances between degenerated $U$-statistics of the same order. The problem of writing the explicit Hoeffding-ANOVA decomposition of a given random variable is usually adressed to characterize the covariance and the consequent asymptotic behavior of such symmetric functionals of the vector $\mathbf{X}_n$, as nondegenerated $U$-statistics or jackknife estimators [see Koroljuk and Borovskich (1994) and Serfling (1980) for a survey], as well as $U$-processes [see, e.g., Arcones and Giné (1993)]. However, it has been completely solved in only two cases: when $\mathbf{X}$ is a sequence of i.i.d. random variables [as first proved in Hoeffding (1948), see, e.g., Hajek (1968), Efron and Stein (1981), Karlin and Rinott (1982), Takemura (1983), Vitale (1990), Bentkus, Götze and van Zwet (1997) and the references therein], and when $\mathbf{X}$ is a collection of $N-1$ extractions without replacement from a finite population [see Zhao and Chen (1990) and Bloznelis and Götze (2001, 2002)], and in both instances, the degenerated $U$-statistics $T_i$ turn out to be linear combinations of well chosen conditional expectations of $T$.

The aim of this paper is twofold.

On the one hand, we shall establish a necessary and sufficient condition for a general exchangeable sequence to be Hoeffding decomposable. Our main result states, indeed, that $\mathbf{X}$ is Hoeffding decomposable if, and only if, $\mathbf{X}$ is composed of *weakly independent* random variables. The notion of weak independence is introduced here for the first time, and will be formally explored in Section 4. To capture the idea of weak independence, suppose $\mathbf{X} = (X_1, X_2, X_3)$, then, $\mathbf{X}$ is weakly independent if, and only if, the following implication holds:

$$\mathbb{E}(\phi(X_1, X_2)|X_1) = 0 \quad \text{a.s.-}\mathbb{P} \quad \Longrightarrow \quad \mathbb{E}(\phi(X_1, X_2)|X_3) = 0 \quad \text{a.s.-}\mathbb{P},$$

where $\phi$ is an arbitrary symmetric kernel such that $\mathbb{E}[\phi(\mathbf{X}_2)^2] < +\infty$. We will see that not every exchangeable sequence is weakly independent, and, therefore, that not every exchangeable sequence is Hoeffding decomposable.

On the other hand, we will apply the above results to explicitly calculate, for every $n$, the Hoeffding-ANOVA decomposition of a general, symmetric $T \in L^2(\mathbf{X}_n)$, when $\mathbf{X}$ is a generalized urn sequence (GUS), a notion that will be introduced in Section 5. As discussed below, the family of GUS contains exclusively exchangeable sequences; examples are i.i.d. random variables, extractions without replacement from a finite population, as well as *generalized Pólya urn schemes* [such as the ones introduced in Ferguson (1973) and Blackwell and MacQueen (1973)]. Consequently, our formulae will extend and unify the classic results about ANOVA decompositions for i.i.d. variables and finite population statistics, and will show that exchangeability is quite a natural framework for studying ANOVA-type decompositions of symmetric statistics. Note, however, that exchangeability is not a necessary condition for a random sequence to be Hoeffding decomposable, see, for example, Karlin and Rinott (1982), Friedrich (1989) and Alberink and Bentkus (1999), where the authors study the case of independent but not identically distributed random variables. In a companion paper [see Peccati (2002a), but also Peccati (2002b, 2003)], we apply our results concerning generalized Pólya urns to obtain a "chaotic decomposition" of the space of square integrable functionals of a Dirichlet–Ferguson process [see, e.g., Ferguson (1973)] defined on a Polish space $(A, \mathcal{A})$.

The paper is organized as follows: in Section 2 we introduce some notation; in Section 3 we define the notion of *Hoeffding spaces* and establish some useful results about exchangeable sequences and (symmetric) $U$-statistics; Section 4 is devoted to the relations between Hoeffding decomposability and weak independence; in Section 5 we prove our main theorems about GUS, whereas Section 6 is devoted to further examples, refinements and applications.

Part of the results of this paper have been announced in Peccati (2003).

## 2. Basic notation.

Fix $n \geq 1$. For any $m \in \{0, 1, \ldots, n\}$, we define

$$V_n(m) := \{\mathbf{k}_{(m)} = (k_1, \ldots, k_m) : 1 \leq k_1 < \cdots < k_m \leq n\}$$

with the convention $\mathbf{k}_{(0)} := 0$ and $V_n(0) = \{0\}$. We also set

$$V_\infty(m) = \bigcup_{n \geq m} V_n(m).$$

For $n \geq m \geq 1$, $\mathbf{l}_{(m)} = (l_1, \ldots, l_m) \in V_\infty(m)$ and $\mathbf{k}_{(n)} = (k_1, \ldots, k_n) \in V_\infty(n)$, $\mathbf{l}_{(m)} \wedge \mathbf{k}_{(n)}$ stands for the class

$$\{l_i : l_i = k_j \text{ for some } j = 1, \ldots, n\}$$

written as an element of $V_\infty(r)$, where $r := \mathrm{Card}\{\mathbf{l}_{(m)} \wedge \mathbf{k}_{(n)}\}$. Analogously, for any $n, m \geq 0$, $\mathbf{k}_{(n)} \setminus \mathbf{l}_{(m)}$ will indicate the set $\{k_j : k_j \neq l_i \ \forall\, i = 1, \ldots, m\}$ written as an element of the class $V_\infty(n-r)$. Again, given $\mathbf{k}_{(n)} \in V_\infty(n)$ and a vector $\mathbf{h}_{(m)} = (h_1, \ldots, h_m)$, by $\mathbf{h}_{(m)} \subset \mathbf{k}_{(n)}$ we will mean that $\mathbf{h}_{(m)} \in V_\infty(m)$, and that for every $i \in \{1, \ldots, m\}$, there exists $j \in \{1, \ldots, n\}$ such that $k_j = h_i$.

As in the Introduction, we now fix $N \in \mathbb{N} \cup \{+\infty\}$ and consider an exchangeable sequence $\mathbf{X} = \{X_n : 1 \leq n < N\}$ composed of random variables with values in the Borel space $(A, \mathcal{A})$. By exchangeability we mean that the law of $\mathbf{X}$ is invariant under finite permutations of the index set $\{n : 1 \leq n < N\}$. More to the point, when $N < +\infty$, $\mathbf{X}$ will always satisfy by convention the following:

ASSUMPTION A. When $N$ is finite, the vector $\mathbf{X} = (X_1, \ldots, X_{N-1})$ is composed of the first $N-1$ elements of a finite exchangeable sequence $(X_1, \ldots, X_{2N-2})$.

In the terminology of Aldous (1983), Assumption A implies that $(X_1, \ldots, X_{N-1})$ is a $2(N-1)$-*extendible* exchangeable sequence. [We recall that, according to Aldous (1983), for $2 \leq M < +\infty$, an exchangeable vector $(Y_1, \ldots, Y_M)$ is said to be $(M+k)$-extendible $(k \geq 1)$ if there exists an exchangeable vector $(Z_1, \ldots, Z_{M+k})$ such that

$$(Y_1, \ldots, Y_M) \overset{\text{law}}{=} (Z_1, \ldots, Z_M).$$

Of course, not every exchangeable vector is extendible.] This point will play an important role in the next section. Recall that if $N = +\infty$, and, therefore, $\mathbf{X}$ is an infinite exchangeable sequence, de Finetti's theorem [see Aldous (1983)] implies that $\mathbf{X}$ is a *mixture* of i.i.d. sequences.

For any $1 \leq n < N$, we define

$$\mathbf{X}_n = (X_1, \ldots, X_n),$$
$$\mathbf{X}_0 = 0$$

and, for any $n \geq 0$ and every $\mathbf{j}_{(n)} \in V_\infty(n)$, we write

$$\mathbf{X}_{\mathbf{j}_{(n)}} = (X_{j_1}, \ldots, X_{j_n}).$$

Now fix $1 \leq n < N$, and consider a symmetric and measurable function $T$ on $A^n$ such that $T(\mathbf{X}_n) \in L^1(\mathbf{X}_n)$. Then, exchangeability implies that for every $0 \leq r \leq m \leq n$, there exists a measurable function

$$[T]_{n,m}^{(r)} : A^m \mapsto \Re$$

with the following properties: (a) for every $\mathbf{j}_{(n)} \in V_\infty(n)$ and $\mathbf{i}_{(m)} \in V_\infty(m)$ satisfying $\mathrm{Card}\{\mathbf{i}_{(m)} \wedge \mathbf{j}_{(n)}\} = r$, one has

$$(1) \qquad \mathbb{E}[T(\mathbf{X}_{\mathbf{j}_{(n)}})|\mathbf{X}_{\mathbf{i}_{(m)}}] = [T]_{n,m}^{(r)}(\mathbf{X}_{\mathbf{i}_{(m)} \wedge \mathbf{j}_{(n)}}, \mathbf{X}_{\mathbf{i}_{(m)} \setminus \mathbf{j}_{(n)}}) \qquad \text{a.s.-}\mathbb{P};$$

(b) for any fixed $(a_1, \ldots, a_r) \in A^r$, the application

$$(a_{r+1}, \ldots, a_m) \mapsto [T]_{n,m}^{(r)}(a_1, \ldots, a_r, a_{r+1}, \ldots, a_m)$$

is symmetric; (c) for any fixed $(a_{r+1}, \ldots, a_m) \in A^{m-r}$, the application

$$(a_1, \ldots, a_r) \mapsto [T]_{n,m}^{(r)}(a_1, \ldots, a_r, a_{r+1}, \ldots, a_m)$$

is symmetric. We will denote by $\widetilde{[T]}_{n,m}^{(r)}$ the canonical symmetrization of $[T]_{n,m}^{(r)}$, that is,

$$
\begin{aligned}
(2) \quad \widetilde{[T]}_{n,m}^{(r)}(a_1, \ldots, a_m) &= \frac{1}{m!} \sum_{\pi} [T]_{n,m}^{(r)}(a_{\pi(1)}, \ldots, a_{\pi(m)}) \\
&= \binom{m}{r}^{-1} \sum_{\mathbf{j}_{(r)} \in V_m(r)} [T]_{n,m}^{(r)}(\mathbf{a}_{\mathbf{j}_{(r)}}, \mathbf{a}_{(1,\ldots,m) \setminus \mathbf{j}_{(r)}}),
\end{aligned}
$$

where $\mathbf{a}_{\mathbf{j}_{(r)}} = (a_{j_1}, \ldots, a_{j_r})$, for every $r \le m$ and every $\mathbf{j}_{(r)} \in V_m(r)$, and $\pi$ runs over all permutations of the set $(1, \ldots, m)$.

## 3. Hoeffding spaces associated to exchangeable sequences.

3.1. *Hoeffding spaces.* Let the previous notation and assumptions prevail throughout this section. For a certain $1 \le n < N$, we introduce the following notation. Set $U_0 = \Re$ and, for $i = 1, \ldots, n$,

$$U_i(\mathbf{X}_n) = \overline{\text{v.s. } \{T(\mathbf{X}_{\mathbf{j}_{(i)}}) : T(\mathbf{X}_{\mathbf{j}_{(i)}}) \in L^2(\mathbf{X}_n), \mathbf{j}_{(i)} \in V_n(i)\}}^{L^2(\mathbf{X}_n)},$$

where v.s.$\{B\}$ indicates the vector space generated by $B$, and eventually

$$H_0 = U_0,$$

$$H_i(\mathbf{X}_n) = U_i(\mathbf{X}_n) \cap U_{i-1}(\mathbf{X}_n)^{\perp}, \qquad i = 1, \ldots, n,$$

where $U_{i-1}(\mathbf{X}_n)^{\perp}$ denotes, for every $i$, the orthogonal of $U_{i-1}(\mathbf{X}_n)$ in $L^2(\mathbf{X}_n)$ [the reader is referred, e.g., to Dudley (1989) for any unexplained notion concerning Hilbert spaces]. We also set $L_s^2(\mathbf{X}_n)$ to be the subspace of $L^2(\mathbf{X}_n)$ composed of symmetric functionals of the vector $\mathbf{X}_n$ and eventually, for $i = 1, \ldots, n$,

$$SU_0 = SH_0 = \Re,$$

$$SU_i(\mathbf{X}_n) = \overline{\text{v.s.} \left\{ T : T = \sum_{\mathbf{j}_{(i)} \in V_n(i)} \phi(\mathbf{X}_{\mathbf{j}_{(i)}}), \phi(\mathbf{X}_i) \in L_s^2(\mathbf{X}_i) \right\}}^{L_s^2(\mathbf{X}_n)},$$

$$SH_i(\mathbf{X}_n) = SU_i(\mathbf{X}_n) \cap SU_{i-1}(\mathbf{X}_n)^{\perp},$$

where, in the last formula, the orthogonal is taken in $L_s^2(\mathbf{X}_n)$.

We define $\{H_i(\mathbf{X}_n): i = 1, \ldots, n\}$ and $\{SH_i(\mathbf{X}_n): i = 1, \ldots, n\}$ to be respectively the collection of *Hoeffding spaces* and *symmetric Hoeffding spaces* associated to $\mathbf{X}_n$. It is immediate that the class $U_i(\mathbf{X}_n)$ represents, for a fixed $i \leq n$, the span of those functionals of $\mathbf{X}_n$ that depend *at most* on $i$ components of the vector $\mathbf{X}_n$, and that the $H_i(\mathbf{X}_n)$'s are obtained as a Gram–Schmidt orthogonalization [see Dudley (1989)] of the increasing sequence $\{U_i(\mathbf{X}_n)\}$. On the other hand, $SU_i(\mathbf{X}_n)$ is the subspace of $U_i(\mathbf{X}_n)$ generated by $U$-statistics, based on $\mathbf{X}_n$, with symmetric and square integrable kernels of order $i$.

Given $T \in L^2(\mathbf{X}_n)$, for every $i = 0, \ldots, n$, we will use the symbols

$$\pi[T, H_i](\mathbf{X}_n) \quad \text{and} \quad \pi[T, SH_i](\mathbf{X}_n)$$

to indicate the projection of $T$ on $H_i(\mathbf{X}_n)$ and $SH_i(\mathbf{X}_n)$. Of course, for every $T \in L^2(\mathbf{X}_n)$,

$$T = \mathbb{E}(T) + \sum_{i=1}^n \pi[T, H_i](\mathbf{X}_n)$$

and for every $T \in L_s^2(\mathbf{X}_n)$,

$$T = \mathbb{E}(T) + \sum_{i=1}^n \pi[T, SH_i](\mathbf{X}_n).$$

The rest of the paper is essentially devoted to the characterization of the operators

$$\pi[\cdot, SH_i](\mathbf{X}_n): L_s^2(\mathbf{X}_n) \mapsto SH_i(\mathbf{X}_n): T \mapsto \pi[T, SH_i](\mathbf{X}_n)$$

for $\mathbf{X}$ belonging to some special class of exchangeable sequences. In particular, we will be interested in sequences satisfying the following:

DEFINITION 1. The exchangeable sequence $\mathbf{X}$ is said to be *Hoeffding decomposable* if, for every $1 \leq n < N$ and every $1 \leq i \leq n$, the following double implication holds: $T \in SH_i(\mathbf{X}_n)$ if, and only if, there exists

$$\phi_T^{(i)}: A^i \mapsto \Re$$

such that $\phi_T^{(i)}(\mathbf{X}_i) \in L_s^2(\mathbf{X}_i)$,

$$(3) \qquad \mathbb{E}[\phi_T^{(i)}(\mathbf{X}_i)|\mathbf{X}_{i-1}] = 0, \qquad \mathbb{P}\text{-a.s.}$$

and

$$T = \sum_{\mathbf{j}_{(i)} \in V_n(i)} \phi_T^{(i)}(\mathbf{X}_{\mathbf{j}_{(i)}}).$$

Of course, the crucial point in the above definition is given by (3). When $\phi_T^{(i)}$ is such that $\phi_T^{(i)}(\mathbf{X}_i) \in L_s^2(\mathbf{X}_i)$ and satisfies (3), we write

$$\phi_T^{(i)} \in \Xi_i(\mathbf{X}).$$

It is well known that i.i.d. sequences are Hoeffding decomposable. As already pointed out, this feature has been the key tool to study the asymptotic behavior of symmetric $U$-statistics, *via* the characterization of their covariance structure [see, e.g., Serfling (1980) and Vitale (1990)]. We will see in the next section that another archetypal class of Hoeffding decomposable sequences is given by extractions without replacement from finite populations.

3.2. *Hoeffding decompositions for finite population statistics.* In this section we shall shortly recall some of the findings of Zhao and Chen (1990) that will be useful in the following sections. Note that the theory of Hoeffding decompositions for finite population statistics has been further developed in the works of Bloznelis and Götze (2001, 2002) that have inspired our presentation.

Fix $M \geq 1$. We note $\mathbf{z} = (z_1, \ldots, z_M)$, a *nonordered collection* of $M$ elements of $A$, and we identify $\mathbf{z}$ with the measure on $(A, \mathcal{A})$ given by

$$\mu_{\mathbf{z}}(C) = \frac{1}{M} \sum_{i=1}^M \mathbf{1}_C(z_i), \qquad C \in \mathcal{A}.$$

We note $\mathcal{Z}_M(A)$, the set of all such $\mathbf{z}$. To each $\mathbf{z} \in \mathcal{Z}_M(A)$, we associate the random vector

$$\mathbf{Y}^{\mu_{\mathbf{z}}} = (Y_1^{\mu_{\mathbf{z}}}, \ldots, Y_M^{\mu_{\mathbf{z}}}) = (z_{\pi^*(1)}, \ldots, z_{\pi^*(M)}),$$

where $\pi^*$ indicates a random permutation, uniformly distributed over all permutations of $(1, \ldots, M)$. In other words, $\mathbf{Y}^{\mu_{\mathbf{z}}}$ has the law of a vector of $M$ extractions without replacement from a finite population whose composition is given by the measure $\mu_{\mathbf{z}}$. The following result, that is essentially due to Zhao and Chen (1990), characterizes the class of symmetric Hoeffding spaces associated to $\mathbf{Y}_m^{\mu_{\mathbf{z}}} = (Y_1^{\mu_{\mathbf{z}}}, \ldots, Y_m^{\mu_{\mathbf{z}}})$, when $m < M$. Of course, $\mathbf{Y}_m^{\mu_{\mathbf{z}}}$ has the law of the first $m$ extractions without replacement from $\mathbf{z}$.

PROPOSITION 1. *Let $T \in L_s^2(\mathbf{Y}_m^{\mu_{\mathbf{z}}})$, where $\mathbf{z} \in \mathcal{Z}_M(A)$ and $m < M$. Then, there exists a unique class of functions*

$$g_{T, \mu_{\mathbf{z}}}^{(i)} : A^i \mapsto \Re, \qquad i = 1, \ldots, m,$$

*that verify for every $i = 1, \ldots, m$,*

$$\mathbb{E}[g_{T, \mu_{\mathbf{z}}}^{(i)}(\mathbf{Y}_{\mathbf{j}_{(i)}}^{\mu_{\mathbf{z}}}) | \mathbf{Y}_{\mathbf{j}_{(i-1)}}^{\mu_{\mathbf{z}}}] = 0$$

*for every* $\mathbf{j}_{(i)} \in V_m(i)$ *and every* $\mathbf{j}_{(i-1)} \in V_m(i-1)$, *and*

$$(4) \qquad \pi[T, SH_i](\mathbf{Y}_m^{\mu_\mathbf{z}}) = \sum_{\mathbf{j}_{(i)} \in V_m(i)} g_{T,\mu_\mathbf{z}}^{(i)}(\mathbf{Y}_{\mathbf{j}_{(i)}}^{\mu_\mathbf{z}}).$$

*Moreover,* $g_{T,\mu_\mathbf{z}}^{(i)} = 0$ *when* $i > M - m$, *and also,*

$$(5) \qquad \mathbb{E}[\pi[T, SH_i](\mathbf{Y}_m^{\mu_\mathbf{z}})^2] = \frac{\binom{m}{i}\binom{M-m}{i}}{\binom{M-i}{i}} \mathbf{1}_{(M-m \geq i)} \mathbb{E}[g_{T,\mu_\mathbf{z}}^{(i)}(\mathbf{Y}_{\mathbf{j}_{(i)}}^{\mu_\mathbf{z}})^2].$$

Formula (4) implies that $\mathbf{Y}_m^{\mu_\mathbf{z}}$ is Hoeffding decomposable. Proposition 1 will be used in the proof of the main result of the following section. The main reason of its usefulness is nested in the following basic result, whose proof can be found, for example, in Aldous (1983).

PROPOSITION 2. *Under the previous notation, let* $\mathbf{X}_M = (X_1, \ldots, X_M)$, $M < +\infty$, *be a finite exchangeable sequence with values in* $(A, \mathcal{A})$. *Then, conditioned on* $\mu_{\mathbf{X}_M}$, *the law of* $\mathbf{X}_M$ *coincides a.s. with that of* $\mathbf{Y}^{\mu_{\mathbf{X}_M}}$, *that is, a.s.-*$\mathbb{P}$, *for every* $C \in \mathcal{A}^{\otimes M}$,

$$\mathbb{P}[\mathbf{X}_M \in C | \mu_{\mathbf{X}_M}] = \frac{1}{M!} \sum_\pi \mathbf{1}_C(X_{\pi(1)}, \ldots, X_{\pi(M)}),$$

*where* $\pi$ *runs over all permutation of* $(1, \ldots, M)$.

3.3. *Representation of U-statistics for exchangeable observations.* To avoid trivialities, from now on we will systematically work under the following:

ASSUMPTION B. *For every* $1 \leq i \leq n < N$, $H_i(\mathbf{X}_n) \neq \{0\}$ *and* $SH_i(\mathbf{X}_n) \neq \{0\}$.

Assumption B excludes, for instance, the case $X_n = X_1$ for every $n \geq 1$. Note that, under Assumption B, for each $1 \leq i < n$ (as usual, given a collection $\{A, A_j : j = 0, 1, \ldots\}$ of Hilbert spaces, we write $A = \bigoplus A_j$ to mean that $A_j \subset A$ for every $j$, $A_j \perp A_i$ for $i \neq j$ and that every $x \in A$ admits the (unique) representation $x = \sum \pi[x, A_j]$, where $\pi$ stands again for the projection operator),

$$U_i(\mathbf{X}_n) = \bigoplus_{a \leq i} H_a(\mathbf{X}_n) \subsetneq L^2(\mathbf{X}_n) = U_n(\mathbf{X}_n) = \bigoplus_{a=0}^n H_a(\mathbf{X}_n),$$

$$SU_i(\mathbf{X}_n) = \bigoplus_{a \leq i} SH_a(\mathbf{X}_n) \subsetneq L_s^2(\mathbf{X}_n) = SU_n(\mathbf{X}_n) = \bigoplus_{a=0}^n SH_a(\mathbf{X}_n).$$

We shall now show that the elements of $SU_i(\mathbf{X}_n)$ have a unique representation. Our key tool will be the following result.

LEMMA 3. *Let $\mathbf{X} = \{X_n : 1 \le n < N\}$ be an exchangeable sequence, satisfying Assumption* A *in the case of a finite $N$, as well as Assumption* B. *Then, there exist constants $k(N, n, i) \in (0, +\infty)$, $1 \le i \le n < N$, depending uniquely on $N$, $n$ and $i$ (and not on the law of $\mathbf{X}$) satisfying for every $i = 1, \ldots, n$, and every real valued $\phi(\cdot)$, defined on $A^i$ and such that $\phi(\mathbf{X}_i) \in L_s^2(\mathbf{X}_i)$,*

$$\mathbb{E}\bigg[\bigg(\sum_{\mathbf{j}_{(i)} \in V_n(i)} \phi(\mathbf{X}_{\mathbf{j}_{(i)}})\bigg)^2\bigg] \ge k(N, n, i)\mathbb{E}[\phi(\mathbf{X}_i)^2].$$

PROOF. We start with the case $N = +\infty$. In this case, de Finetti's theorem [see once again Aldous (1983), Section 7] yields the existence of a random probability measure $D(\cdot; \omega)$ such that, conditioned to $D$, $\mathbf{X}$ is a sequence of i.i.d. random variables with common law equal to $D$. It follows that [noting $\binom{a}{b}_* = \binom{a}{b}\mathbf{1}_{(a \ge b)}$], due to symmetry and exchangeability,

$$\mathbb{E}\bigg[\bigg(\sum_{\mathbf{j}_{(i)} \in V_n(i)} \phi(\mathbf{X}_{\mathbf{j}_{(i)}})\bigg)^2\bigg]$$

$$= \sum_{\mathbf{j}_{(i)} \in V_n(i)} \sum_{r=0}^{i} \binom{i}{r}\binom{n-i}{i-r}_* \mathbb{E}[\phi(\mathbf{X}_i)\phi(\mathbf{X}_r, X_{i+1}, \ldots, X_{2i-r})]$$

$$= \sum_{\mathbf{j}_{(i)} \in V_n(i)} \sum_{r=0}^{i} \binom{i}{r}\binom{n-i}{i-r}_* \mathbb{E}\bigg[\int_{A^r} D^{\otimes r}(da_1, \ldots, da_r)$$

$$\times \bigg(\int_{A^{i-r}} D^{\otimes i-r}(da_{i+1}, \ldots, da_{2i-r})$$

$$\times \phi(a_1, \ldots, a_r, a_{i+1}, \ldots, a_{2i-r})\bigg)^2\bigg]$$

$$\ge \sum_{\mathbf{j}_{(i)} \in V_n(i)} \mathbb{E}[\phi(\mathbf{X}_{\mathbf{j}_{(i)}})^2]$$

$$= \binom{n}{i}\mathbb{E}[\phi(\mathbf{X}_i)^2].$$

Now we deal with a finite $N$. We recall that, in our setting, $\mathbf{X}$ is in this case of the form $(X_1, \ldots, X_{N-1})$, with $\mathbf{X}_{2N-2} = (X_1, \ldots, X_{2N-2})$ indicating an exchangeable vector of $2(N-1)$ random variables. Then, we use extensively the content and the notation of Propositions 1 and 2 to obtain, due again to symmetry and exchangeability,

$$\mathbb{E}\bigg[\bigg(\sum_{\mathbf{j}_{(i)} \in V_n(i)} \phi(\mathbf{X}_{\mathbf{j}_{(i)}})\bigg)^2\bigg]$$

$$= \mathbb{E}\left[\mathbb{E}\left[\left(\sum_{\mathbf{j}_{(i)} \in V_n(i)} \phi(\mathbf{X}_{\mathbf{j}_{(i)}})\right)^2 \Big| \mu_{\mathbf{X}_{2N-2}}\right]\right]$$

(6)
$$= \mathbb{E}\left[\mathbb{E}\left[\left(\binom{n}{i} \mathbb{E}[\phi(\mathbf{X}_i)|\mu_{\mathbf{X}_{2N-2}}]\right.\right.\right.$$
$$\left.\left.\left. + \sum_{k=1}^{i} \sum_{\mathbf{j}_{(k)} \in V_n(k)} \binom{n-k}{i-k} g^{(k)}_{\phi,\mu_{\mathbf{X}_{2N-2}}}(\mathbf{X}_{\mathbf{j}_{(k)}})\right)^2 \Big| \mu_{\mathbf{X}_{2N-2}}\right]\right]$$

and, therefore,

$$\mathbb{E}\left[\left(\sum_{\mathbf{j}_{(i)} \in V_n(i)} \phi(\mathbf{X}_{\mathbf{j}_{(i)}})\right)^2\right]$$

(7)
$$= \mathbb{E}\left[\binom{n}{i}^2 \mathbb{E}[\phi(\mathbf{X}_i)|\mu_{\mathbf{X}_{2N-2}}]^2\right.$$
$$\left. + \sum_{k=1}^{i} \frac{\binom{n-k}{i-k}^2 \binom{2N-2-n}{k}\binom{n}{k}}{\binom{2N-2-k}{k}} \mathbb{E}[g^{(k)}_{\phi,\mu_{\mathbf{X}_{2N-2}}}(\mathbf{X}_k)^2|\mu_{\mathbf{X}_{2N-2}}]\right].$$

To be clear, the calculations contained in (6) and (7) are performed as follows. First, write the Hoeffding decomposition of $\phi(\mathbf{X}_{\mathbf{j}_{(i)}})$, under the conditioned probability $\mathbb{P}[\cdot|\mu_{\mathbf{X}_{2N-2}}]$ and for every $\mathbf{j}_{(i)} \in V_n(i)$. Then, by using the relation

$$\mathbb{E}[g^{(k)}_{\phi,\mu_{\mathbf{X}_{2N-2}}}(\mathbf{X}_{\mathbf{j}_{(k)}})|\mu_{\mathbf{X}_{2N-2}}, \mathbf{X}_{\mathbf{j}_{(k-1)}}] = 0, \qquad \mathbb{P}[\cdot|\mu_{\mathbf{X}_{2N-2}}]\text{-a.s.}$$

for every $\mathbf{j}_{(k-1)} \in V_n(k-1)$ [that can be verified directly, by inspecting the proof of the main results of Bloznelis and Götze (2001) or by using Corollary 9; i.e., not circular reasoning, as a matter of fact, to prove Proposition 8 and Corollary 9, we do not need Lemma 3], observe that

$$\sum_{\mathbf{j}_{(k)} \in V_n(k)} \binom{n-k}{i-k} g^{(k)}_{\phi,\mu_{\mathbf{X}_{2N-2}}}(\mathbf{X}_{\mathbf{j}_{(k)}})$$

is the projection of $\sum_{\mathbf{j}_{(i)} \in V_n(i)} \phi(\mathbf{X}_{\mathbf{j}_{(i)}})$ on the $k$th symmetric Hoeffding space associated to $\mathbf{X}_n$ under the measure $\mathbb{P}[\cdot|\mu_{\mathbf{X}_{2N-2}}]$. Finally, use Proposition 1. Now write

$$k(N, n, i) = \min\left\{\binom{n}{i}^2; \frac{\binom{n-s}{i-s}^2\binom{2N-2-n}{s}\binom{n}{s}}{\binom{2N-2-s}{s}\binom{i}{s}^2}, s = 1, \ldots, i\right\}$$

to obtain, thanks to the Jensen inequality,

$$\mathbb{E}\left[\left(\sum_{\mathbf{j}_{(i)}\in V_n(i)}\phi(\mathbf{X}_{\mathbf{j}_{(i)}})\right)^2\right]$$

$$\geq k(N,n,i)\mathbb{E}\left[\mathbb{E}[\phi(\mathbf{X}_i)|\mu_{\mathbf{X}_{2N-2}}]^2\right.$$

$$\left.+\sum_{k=1}^{i}\binom{i}{k}^2\mathbb{E}[g_{\phi,\mu_{\mathbf{X}_{2N-2}}}^{(k)}(\mathbf{X}_k)^2|\mu_{\mathbf{X}_{2N-2}}]\right]$$

$$\geq k(N,n,i)\mathbb{E}\left[\mathbb{E}[\phi(\mathbf{X}_i)|\mu_{\mathbf{X}_{2N-2}}]^2\right.$$

$$\left.+\sum_{k=1}^{i}\mathbb{E}\left[\left(\sum_{\mathbf{j}_{(k)}\in V_i(k)}g_{\phi,\mu_{\mathbf{X}_{2N-2}}}^{(k)}(\mathbf{X}_{\mathbf{j}_{(k)}})\right)^2\Big|\mu_{\mathbf{X}_{2N-2}}\right]\right]$$

$$=k(N,n,i)\mathbb{E}[\mathbb{E}[\phi(\mathbf{X}_i)^2|\mu_{\mathbf{X}_{2N-2}}]],$$

which yields the desired result. $\square$

REMARK. An inspection of the proof of Lemma 3 shows the relevance of the assumption: for a finite $N$, $\mathbf{X}=(X_1,\ldots,X_{N-1})$ is a $2(N-1)$-extendible sequence. Suppose indeed that $(X_1,\ldots,X_{N-1})$ are the first $N-1$ instants of a sequence $\mathbf{X}_M=(X_1,\ldots,X_M)$, with $N\leq M<2N-2$. Then, according to Proposition 1,

$$\sum_{k=1}^{i}\sum_{\mathbf{j}_{(k)}\in V_n(k)}\binom{n-k}{i-k}g_{\phi,\mu_{\mathbf{X}_M}}^{(k)}(\mathbf{X}_{\mathbf{j}_{(k)}})$$

$$=\sum_{k=1}^{\min\{i,M-i\}}\sum_{\mathbf{j}_{(k)}\in V_n(k)}\binom{n-k}{i-k}g_{\phi,\mu_{\mathbf{X}_M}}^{(k)}(\mathbf{X}_{\mathbf{j}_{(k)}}),\qquad\text{a.s.-}\mathbb{P}(\cdot|\mu_{\mathbf{X}_M}),$$

and, also,

$$\mathbb{E}\left[\sum_{k=1}^{\min\{i,M-i\}}\left(\sum_{\mathbf{j}_{(k)}\in V_n(k)}\binom{n-k}{i-k}g_{\phi_{\mathbf{X}_M}}^{(k)}(\mathbf{X}_{\mathbf{j}_{(k)}})\right)^2\Big|\mu_{\mathbf{X}_M}\right]$$

$$\tag{8}$$

$$=\sum_{k=1}^{\min\{i,M-n\}}\frac{\binom{n-k}{i-k}^2\binom{M-n}{k}\binom{n}{k}}{\binom{M-k}{k}}\mathbb{E}[g_{\phi,\mu_{\mathbf{X}_M}}^{(k)}(\mathbf{X}_k)^2|\mu_{\mathbf{X}_M}].$$

It is easily seen that, when $i\geq M-i>M-n$, relation (8) does not allow to conclude the proof of Lemma 3.

Lemma 3 has important consequences which are stated in the next two corollaries.

COROLLARY 4.   *For $1 \le i \le n < N$, suppose the applications $\phi$ and $\phi'$, both from $A^i$ to $\Re$, are such that $\phi(\mathbf{X}_i), \phi'(\mathbf{X}_i) \in L^2_s(\mathbf{X}_i)$. Then,*

$$\sum_{\mathbf{j}_{(i)} \in V_n(i)} \phi(\mathbf{X}_{\mathbf{j}_{(i)}}) = \sum_{\mathbf{j}_{(i)} \in V_n(i)} \phi'(\mathbf{X}_{\mathbf{j}_{(i)}}), \qquad \mathbb{P}\text{-}a.s.$$

*implies*

$$\phi(\mathbf{X}_i) = \phi'(\mathbf{X}_i), \qquad \mathbb{P}\text{-}a.s.$$

Corollary 4 says that elements of $SU_i(\mathbf{X}_n)$ admit an essentially unique representation as $U$-statistics with symmetric kernel of order $i$. The next result states that $SU_i(\mathbf{X}_n)$ contains exclusively random variables of this kind.

COROLLARY 5.   *For $1 \le i \le n < N$,*

$$SU_i(\mathbf{X}_n) = \left\{ T : T = \sum_{\mathbf{j}_{(i)} \in V_n(i)} \phi(\mathbf{X}_{\mathbf{j}_{(i)}}), \phi(\mathbf{X}_i) \in L^2_s(\mathbf{X}_i) \right\}.$$

PROOF.   For fixed $i$ and $n$ as in the statement, just observe that if the family $\{T^{(l)} : l \ge 1\}$, defined as

$$T^{(l)} = \sum_{\mathbf{j}_{(i)} \in V_n(i)} \phi^{(l)}(\mathbf{X}_{\mathbf{j}_{(i)}}), \qquad \phi^{(l)}(\mathbf{X}_i) \in L^2_s(\mathbf{X}_i), l = 1, 2, \ldots,$$

is a Cauchy sequence in $L^2_s(\mathbf{X}_n)$, then Lemma 3 implies that $\phi^{(l)}(\mathbf{X}_i)$ is also Cauchy in $L^2_s(\mathbf{X}_i)$.   □

**4. Hoeffding decomposability and weak independence.**   For the rest of the section $\mathbf{X}$ will be a possibly infinite exchangeable sequence satisfying both Assumptions A and B.

DEFINITION 2.   We say that the sequence $\mathbf{X}$ is composed of *weakly independent* random variables (or that the sequence $\mathbf{X}$ is *weakly independent*) if for every $1 \le n < N$ and every $T \in L^2_s(\mathbf{X}_n)$,

$$[T]^{(n-1)}_{n,n-1}(\mathbf{X}_{n-1}) = 0, \qquad \text{a.s.-}\mathbb{P},$$

*implies*

$$\widetilde{[T]}^{(r)}_{n,n-1}(\mathbf{X}_{n-1}) = 0, \qquad \text{a.s.-}\mathbb{P},$$

for every $0 \le r \le n-1$ such that $2n - r \le N$, where the functions $[T]^{(\cdot)}_{\cdot,\cdot}$ and $\widetilde{[T]}^{(\cdot)}_{\cdot,\cdot}$ have been introduced, respectively, in (1) and (2).

Of course, independence implies weak independence. Another example of weak independence is given by sampling without replacement and, in general, by the class of GUS that we will discuss in the next section. However, not every exchangeable sequence is weakly independent.

EXAMPLE (A class of exchangeable sequences that are not weakly independent). Consider an infinite sequence

$$\mathbf{X} = \{X_n : n \geq 1\}$$

with values in $\{0, 1\}$, whose law is determined by the following relation, valid for every $n \geq 1$ and every $(e_1, \ldots, e_n) \in \{0, 1\}^n$:

$$\mathbb{P}(X_1 = e_1, \ldots, X_n = e_n) = \varepsilon^{-1} \int_0^\varepsilon \prod_{i=1}^n x^{e_i} (1-x)^{1-e_i} \, dx,$$

where $\varepsilon$ is a fixed constant such that $0 < \varepsilon < 1$.

This is equivalent to saying that, conditioned on the realization of a real valued random variable $Y$ such that

$$\mathbb{P}(Y \in C) = \varepsilon^{-1} \int_{(0,\varepsilon) \cap C} dy,$$

the sequence $\mathbf{X}$ is composed of independent Bernoulli trials with common parameter equal to $Y$. In this case, a necessary condition for $\mathbf{X}$ to be weakly independent is that for any symmetric $\phi$ on $\{0, 1\}^2$ such that

(9) $$\mathbb{E}(\phi(X_1, X_2) | X_2) = 0$$

must also hold

(10) $$\mathbb{E}(\phi(X_1, X_2) | X_3) = 0.$$

We shall construct a symmetric $\phi$ that respects (9) but not (10). Define, indeed,

$$\phi(1, 0) = \phi(0, 1) = 1$$

and also

$$\phi(1, 1) = -\frac{\int_0^\varepsilon x(1-x) \, dx}{\int_0^\varepsilon x^2 \, dx} = 1 - \frac{3}{2\varepsilon},$$

$$\phi(0, 0) = -\frac{\int_0^\varepsilon x(1-x) \, dx}{\int_0^\varepsilon (1-x)^2 \, dx} = \frac{\varepsilon^2 - (3/2)\varepsilon}{3 - 3\varepsilon + \varepsilon^2},$$

so that

$$\mathbb{E}(\phi(X_1, X_2) | X_2 = 0) = \mathbb{E}(\phi(X_1, X_2) | X_2 = 1) = 0,$$

and also

$$\mathbb{E}(\phi(X_1, X_2)|X_3 = 0) = \frac{1}{8} \frac{\varepsilon^3(\varepsilon - 1)}{(3 - 3\varepsilon + \varepsilon^2)(\varepsilon - (1/2)\varepsilon^2)} < 0,$$

since $\varepsilon \in (0, 1)$.

It is interesting to note that by taking $\varepsilon = 1$, one would obtain a weakly independent sequence. As a matter of fact, $\mathbf{X}$ is in this case a Pólya urn sequence with parameters $(1, 1)$ (see the discussion below).

The following result establishes a necessary and sufficient condition for Hoeffding decomposability.

THEOREM 6. *The exchangeable sequence $\mathbf{X}$ is Hoeffding decomposable if, and only if, it is weakly independent.*

PROOF. To simplify, we will systematically consider r.v.'s $T$ such that $\mathbb{E}(T) = 0$. Now suppose that the sequence $\mathbf{X}$ is weakly independent, and take $T(\mathbf{X}_2) \in L^2_s(\mathbf{X}_2)$. According to Corollary 5, there exists a function $\phi_T^{(1)} : A \mapsto \Re$ such that $\mathbb{E}[(\phi_T^{(1)}(X_1))^2] < +\infty$, and also

$$\pi[T, SH_1](\mathbf{X}_2) = \phi_T^{(1)}(X_1) + \phi_T^{(1)}(X_2),$$

(11)      $$\pi[T, SH_2](\mathbf{X}_2) = T(\mathbf{X}_2) - \phi_T^{(1)}(X_1) - \phi_T^{(1)}(X_2)$$

$$= \phi_T^{(2)}(\mathbf{X}_2).$$

Plainly, $\phi_T^{(2)} \in \Xi_2(\mathbf{X})$: as a matter of fact, for every bounded $h$ on $A$ and thanks to exchangeability and symmetry,

$$\mathbb{E}[\phi_T^{(2)}(\mathbf{X}_2)h(X_1)] = \tfrac{1}{2}\mathbb{E}[\phi_T^{(2)}(\mathbf{X}_2)(h(X_1) + h(X_2))]$$

$$= 0.$$

Now take $n > 2$. To show that if $G \in SH_2(\mathbf{X}_n)$, then there exists $\phi_G^{(2)} \in \Xi_2(\mathbf{X})$ such that

(12)      $$G = \sum_{\mathbf{j}_{(2)} \in V_n(2)} \phi_G^{(2)}(\mathbf{X}_{\mathbf{j}_{(2)}}),$$

it is sufficient to show that representation (12) holds for random variables of the type

$$G = \pi[F, SH_2](\mathbf{X}_n),$$

where $F$ is centered and such that $F \in SU_2(\mathbf{X}_n)$. Thanks again to Corollary 5, we know that there exists a symmetric and square integrable kernel $T$

such that

$$F = \sum_{\mathbf{j}_{(2)} \in V_n(2)} T(\mathbf{X}_{\mathbf{j}_{(2)}})$$

and also, with the notation introduced in (11),

$$\pi[F, SH_1](\mathbf{X}_n) = (n-1) \sum_{i=1}^{n} \phi_T^{(1)}(X_i),$$

$$\pi[F, SH_2](\mathbf{X}_n) = \sum_{\mathbf{j}_{(2)} \in V_n(2)} \phi_T^{(2)}(\mathbf{X}_{\mathbf{j}_{(2)}}).$$

As a matter of fact,

$$\begin{aligned}
F &= \sum_{\mathbf{j}_{(2)} \in V_n(2)} [T(\mathbf{X}_{\mathbf{j}_{(2)}}) - \phi_T^{(1)}(X_{j_1}) - \phi_T^{(1)}(X_{j_2})] \\
&\quad + \sum_{\mathbf{j}_{(2)} \in V_n(2)} [\phi_T^{(1)}(X_{j_1}) + \phi_T^{(1)}(X_{j_2})] \\
&= \sum_{\mathbf{j}_{(2)} \in V_n(2)} \phi_T^{(2)}(\mathbf{X}_{\mathbf{j}_{(2)}}) + (n-1) \sum_{i=1}^{n} \phi_T^{(1)}(X_i).
\end{aligned}$$

Moreover, for every $h$ such that $\mathbb{E}(h(X_1)^2) < +\infty$,

$$\mathbb{E}\left( \sum_{\mathbf{j}_{(2)} \in V_n(2)} \phi_T^{(2)}(\mathbf{X}_{\mathbf{j}_{(2)}}) \sum_{i=1}^{n} h(X_i) \right) = 0,$$

since we have assumed that $\mathbf{X}$ is weakly independent. Now we use a recurrence argument. Suppose, indeed, that there exists $k \geq 1$ with the following property: for every $k \leq n < N$, for $i = 1, \ldots, k-1$, $F \in SH_i(\mathbf{X}_n)$ implies that there exists $\phi_F^{(i)} \in \Xi_i(\mathbf{X})$ such that

$$(13) \qquad F = \sum_{\mathbf{j}_{(i)} \in V_n(i)} \phi_F^{(i)}(\mathbf{X}_{\mathbf{j}_{(i)}}),$$

and observe that we have verified such a claim for $k = 1, 2, 3$. Given $k$, we shall verify that for every $n \geq k$, a random variable of the type

$$G = \pi[F, SH_k](\mathbf{X}_n)$$

for a generic $F \in SU_k(\mathbf{X}_n)$ has the representation (13) for $i = k$ and $\phi_F^{(i)} \in \Xi_i(\mathbf{X})$. To see this, start with $n = k$, and take a symmetric and square integrable kernel $T$ such that $\mathbb{E}(T(\mathbf{X}_n)) = 0$. Then, there exist $\phi_T^{(i)} \in \Xi_i(\mathbf{X})$,

$i = 1, \ldots, k-1$, such that

$$\pi[T, SH_i](\mathbf{X}_k) = \sum_{\mathbf{j}_{(i)} \in V_k(i)} \phi_T^{(i)}(\mathbf{X}_{\mathbf{j}_{(i)}}), \qquad i = 1, \ldots, k-1,$$

$$(14) \qquad \pi[T, SH_k](\mathbf{X}_k) = T(\mathbf{X}_k) - \sum_{i=1}^{k-1} \sum_{\mathbf{j}_{(i)} \in V_k(i)} \phi_T^{(i)}(\mathbf{X}_{\mathbf{j}_{(i)}})$$

$$= \phi_T^{(k)}(\mathbf{X}_k).$$

Since for every bounded and symmetric function $h$ on $A^{k-1}$ with the form

$$h(a_1, \ldots, a_{k-1}) = \sum_{\pi} \prod_{j=1}^{k-1} \mathbf{1}_{C_j}(a_{\pi(j)}),$$

where $C_1, \ldots, C_{k-1} \in \mathcal{A}$ and $\pi$ runs over all permutations of $(1, \ldots, k-1)$, we have

$$0 = \mathbb{E}\left[\phi_T^{(k)}(\mathbf{X}_k) \sum_{\mathbf{j}_{(k-1)} \in V_k(k-1)} h(\mathbf{X}_{\mathbf{j}_{(k-1)}})\right]$$

$$= k\mathbb{E}[\phi_T^{(k)}(\mathbf{X}_k) h(\mathbf{X}_{k-1})]$$

$$= k!\mathbb{E}\left[\phi_T^{(k)}(\mathbf{X}_k) \prod_{j=1}^{k-1} \mathbf{1}_{C_j}(X_j)\right],$$

due to exchangeability and to the symmetry of $\phi_T^{(k)}$, we obtain immediately $\phi_T^{(k)} \in \Xi_k(\mathbf{X})$. Now, for $n > k$, take $F \in SU_k(\mathbf{X}_n)$ with the form

$$F = \sum_{\mathbf{j}_{(k)} \in V_n(k)} T(\mathbf{X}_{\mathbf{j}_{(k)}}),$$

where $T$ is a centered, square integrable and symmetric kernel. Then, by using the same notation as in (14),

$$F = \sum_{\mathbf{j}_{(k)} \in V_n(k)} \phi_T^{(k)}(\mathbf{X}_{\mathbf{j}_{(k)}}) + \sum_{i=1}^{k-1} \sum_{\mathbf{j}_{(i)} \in V_n(i)} \binom{n-i}{k-i} \phi_T^{(i)}(\mathbf{X}_{\mathbf{j}_{(i)}})$$

and, moreover, for every $h$ on $A^{k-1}$ such that $h(\mathbf{X}_{k-1}) \in L_s^2(\mathbf{X}_{k-1})$,

$$\mathbb{E}\left[\sum_{\mathbf{j}_{(k)} \in V_n(k)} \phi_T^{(k)}(\mathbf{X}_{\mathbf{j}_{(k)}}) \sum_{\mathbf{j}_{(k-1)} \in V_n(k-1)} h(\mathbf{X}_{\mathbf{j}_{(k-1)}})\right] = 0,$$

since for every $\mathbf{j}_{(k-1)} \in V_n(k-1)$,

$$\sum_{\mathbf{j}_{(k)} \in V_n(k)} \mathbb{E}[\phi_T^{(k)}(\mathbf{X}_{\mathbf{j}_{(k)}})|\mathbf{X}_{\mathbf{j}_{(k-1)}}]$$

$$= \sum_{r=0}^{k-1} \sum_{\mathbf{j}_{(k)} \in V_n(k)} \mathbf{1}_{(\mathrm{Card}(\mathbf{j}_{(k)} \wedge \mathbf{j}_{(k-1)})=r)} [\phi_T^{(k)}]_{k,k-1}^{(r)}(\mathbf{X}_{\mathbf{j}_{(k)} \wedge \mathbf{j}_{(k-1)}}, \mathbf{X}_{\mathbf{j}_{(k-1)} \setminus \mathbf{j}_{(k)}})$$

$$= \sum_{r=0}^{k-1} \sum_{\mathbf{j}_{(r)} \subset \mathbf{j}_{(k-1)}} \sum_{\mathbf{j}_{(k)} \in V_n(k)} \mathbf{1}_{(\mathbf{j}_{(k)} \wedge \mathbf{j}_{(k-1)}=\mathbf{j}_{(r)})} [\phi_T^{(k)}]_{k,k-1}^{(r)}(\mathbf{X}_{\mathbf{j}_{(r)}}, \mathbf{X}_{\mathbf{j}_{(k-1)} \setminus \mathbf{j}_{(r)}})$$

and, therefore,

$$\sum_{\mathbf{j}_{(k)} \in V_n(k)} \mathbb{E}[\phi_T^{(k)}(\mathbf{X}_{\mathbf{j}_{(k)}})|\mathbf{X}_{\mathbf{j}_{(k-1)}}]$$

$$= \sum_{r=0}^{k-1} \sum_{\mathbf{j}_{(r)} \subset \mathbf{j}_{(k-1)}} \binom{n-k+1}{k-r}_* [\phi_T^{(k)}]_{k,k-1}^{(r)}(\mathbf{X}_{\mathbf{j}_{(r)}}, \mathbf{X}_{\mathbf{j}_{(k-1)} \setminus \mathbf{j}_{(r)}})$$

$$= \sum_{r=0}^{k-1} \binom{n-k+1}{k-r}_* \binom{k-1}{r} \widetilde{[\phi_T^{(k)}]}_{k,k-1}^{(r)}(\mathbf{X}_{\mathbf{j}_{(k-1)}})$$

$$= 0$$

thanks to the assumption of weak independence and to the fact that $\phi_T^{(k)} \in \Xi_k(\mathbf{X})$.

On the other hand, it is clear that if $\mathbf{X}$ is weakly independent and, for $1 \le i \le n < N$, $F$ has the representation (13) for $\phi^{(i)} \in \Xi_i(\mathbf{X})$, then for any $\mathbf{j}_{(i-1)} \in V_n(i-1)$,

$$\mathbb{E}[F|\mathbf{X}_{\mathbf{j}_{(i-1)}}] = \sum_{r=0}^{i-1} \binom{n-i+1}{i-r}_* \binom{i-1}{r} \widetilde{[\phi_T^{(i)}]}_{i,i-1}^{(r)}(\mathbf{X}_{\mathbf{j}_{(i-1)}}) = 0,$$

and, therefore, $F \in SH_i(\mathbf{X}_n)$.

Thus, we have shown that weak independence implies Hoeffding decomposability. To deal with the opposite implication, suppose for the moment that $N = +\infty$, and that $\mathbf{X}$ is Hoeffding decomposable in the sense of Definition 1. For a given $k \ge 1$, consider a certain $T(\mathbf{X}_k) \in L_s^2(\mathbf{X}_k)$ such that

$$[T]_{k,k-1}^{(k-1)}(\mathbf{X}_{k-1}) = 0, \qquad \mathbb{P}\text{-a.s.}$$

Then,

$$F(\mathbf{X}_{k+1}) = \sum_{\mathbf{j}_{(k)} \in V_{k+1}(k)} T(\mathbf{X}_{\mathbf{j}_{(k)}}) \in SH_k(\mathbf{X}_{k+1})$$

that yields, due to exchangeability and symmetry,

$$0 = \mathbb{E}[F(\mathbf{X}_{k+1})|\mathbf{X}_{k-1}]$$

$$= \sum_{\mathbf{j}_{(k)} \in V_{k+1}(k)} \mathbb{E}[T(\mathbf{X}_{\mathbf{j}_{(k)}})|\mathbf{X}_{k-1}]$$

$$= \sum_{r=k-2} \sum_{\mathbf{j}_{(r)} \in V_{k-1}(r)} \sum_{\mathbf{j}_{(k)} \in V_{k+1}(k)} \mathbf{1}_{(\mathbf{j}_{(k)} \wedge (1,\ldots,k-1) = \mathbf{j}_{(r)})}$$

$$\times [T]_{k,k-1}^{(r)}(\mathbf{X}_{\mathbf{j}_{(r)}}, \mathbf{X}_{(1,\ldots,k-1)\setminus \mathbf{j}_{(r)}})$$

$$= \sum_{\mathbf{j}_{(k-2)} \in V_{k-1}(k-2)} [T]_{k,k-1}^{(k-2)}(\mathbf{X}_{\mathbf{j}_{(k-2)}}, \mathbf{X}_{(1,\ldots,k-1)\setminus \mathbf{j}_{(k-2)}})$$

$$= \frac{(k-1)!}{(k-2)!} [\widetilde{T}]_{k,k-1}^{(k-2)}(\mathbf{X}_{k-1}).$$

Now we use again a recurrence argument. Suppose, indeed, that the Hoeffding decomposability of $\mathbf{X}$ implies the following relation for every $T(\mathbf{X}_k) \in L_s^2(\mathbf{X}_k)$:

$$[T]_{k,k-1}^{(k-1)} = 0 \quad \Longrightarrow \quad [\widetilde{T}]_{k,k-1}^{(k-l)} = 0,$$

for a certain $2 \leq j \leq k-1$, and every $2 \leq l \leq j$. Then, if $T$ is such that $[T]_{k,k-1}^{(k-1)} = 0$, we must have

$$F(\mathbf{X}_{k+j}) = \sum_{\mathbf{j}_{(k)} \in V_{k+j}(k)} T(\mathbf{X}_{\mathbf{j}_{(k)}}) \in SH_k(\mathbf{X}_{k+j}),$$

that implies, again by exchangeability and symmetry,

$$0 = \mathbb{E}[F(\mathbf{X}_{k+j})|\mathbf{X}_{k-1}]$$

$$= \sum_{\mathbf{j}_{(k)} \in V_{k+j}(k)} \mathbb{E}[T(\mathbf{X}_{\mathbf{j}_{(k)}})|\mathbf{X}_{k-1}]$$

$$= \sum_{r=k-j-1} \sum_{\mathbf{j}_{(r)} \in V_{k-1}(r)} \sum_{\mathbf{j}_{(k)} \in V_{k+j}(k)} \mathbf{1}_{(\mathbf{j}_{(k)} \wedge (1,\ldots,k-1) = \mathbf{j}_{(r)})}$$

$$\times [T]_{k,k-1}^{(r)}(\mathbf{X}_{\mathbf{j}_{(r)}}, \mathbf{X}_{(1,\ldots,k-1)\setminus \mathbf{j}_{(r)}})$$

$$= \sum_{r=k-j-1}^{k-1} \binom{j+1}{k-r} \binom{k-1}{r} [\widetilde{T}]_{k,k-1}^{(r)}(\mathbf{X}_{\mathbf{j}_{(r)}}, \mathbf{X}_{(1,\ldots,k-1)\setminus \mathbf{j}_{(r)}})$$

$$= \binom{k-1}{k-j-1} [\widetilde{T}]_{k,k-1}^{(k-j-1)}(\mathbf{X}_{k-1})$$

and, therefore, the desired result. To deal with the case of a finite $N$, just repeat the same argument for $j$ such that $k + j \leq N - 1$. □

One immediate consequence of Theorem 6 is the following:

COROLLARY 7. *Let the exchangeable sequence* $\mathbf{X}$ *be weakly independent. Then, for every* $1 \leq n < N$, *every* $T(\mathbf{X}_n) \in L^2_s(\mathbf{X}_n)$ *and every* $i = 1, \ldots, n$,

$$\pi[T, SH_i](\mathbf{X}_n) = \pi[T, H_i](\mathbf{X}_n).$$

Starting from the next section we analyze the specific case of GUS.

**5. The case of GUS.** In this section we shall investigate the case of GUS, which represent a fundamental example of Hoeffding decomposable sequences. We will consider uniquely the case: $(A, \mathcal{A})$ is a Polish space endowed with its Borel $\sigma$-field. More precisely, for $N \in \mathbb{N} \cup \{\infty\}$, and writing $\mathcal{M}(A)$ for the class of finite and positive measures on $A$, we say that a sequence

$$\mathbf{X}^{(\alpha,c)} = \{X_n^{(\alpha,c)} : 1 \leq n < N\}$$

is a GUS of parameters $\alpha \in \mathcal{M}(A)$ and $c \in \Re$, if $\alpha(A) + c(N-1) > 0$ and if, for every $k$ and every $\mathbf{j}_{(k)} \in V_{N-1}(k)$,

$$(15) \quad \mathbb{P}(X_{j_1}^{(\alpha,c)} \in dx_1, \ldots, X_{j_k}^{(\alpha,c)} \in dx_k) = \prod_{i=1}^{k} \frac{\alpha(dx_i) + \sum_{h=1}^{i-1} \delta_{x_h}^c(dx_i)}{\alpha(A) + c(i-1)},$$

where $\delta_x^c(\cdot) := c\delta_x(\cdot)$, with $\delta_x(\cdot)$ the Dirac measure concentrated in $x$. Note that (15) is equivalent to the following relation: for every $C \in \mathcal{A}$ and for every $n < N$,

$$(16) \quad \mathbb{P}(X_n^{(\alpha,c)} \in C | X_1, \ldots, X_{n-1}) = \frac{\alpha(C) + \sum_{h=1}^{n-1} \delta_{X_h}^c(C)}{\alpha(A) + c(n-1)}.$$

Equations (15) and (16) imply that, for every choice of $\alpha$ and $c$ s.t. $\alpha(A) + c(N-1) > 0$, the sequence $\mathbf{X}^{(\alpha,c)}$ is exchangeable. One can think of $A$ as an urn whose composition is determined by the measure $\alpha(\cdot)$ (thus, $A$ could contain a "continuum" of balls), whereas $\mathbf{X}^{(\alpha,c)}$ represents a sequence of extractions from $A$ according to the following procedure: at each step, one ball is extracted, and $(1 + c)$ balls of the same color are placed in $A$ before the subsequent extraction (one should substitute "placed in" with "eliminated from" when $c < -1$). Note that the assumption $\alpha(A) + c(N-1) > 0$ ensures that *the urn is not exhausted* before the $(N-1)$st step; more to this point: when $c = 0$, $\mathbf{X}^{(\alpha,c)}$ is a sequence of i.i.d. variables with common law $\alpha(\cdot)/\alpha(A)$; if $A = \{a_1, \ldots, a_S\}$, $\alpha$ is the counting measure and $c = -1$, then we must have $\alpha(A) = S > N - 1$ and $\mathbf{X}^{(\alpha,c)}$ has the law of the first $N-1$ extractions without replacement from the finite population $\{a_1, \ldots, a_S\}$ [this is the case studied in Zhao and Chen (1990) and Bloznelis and Götze (2001, 2002)]; when $c > 0$ and $N$ is infinite, $\mathbf{X}^{(\alpha,c)}$ is a *generalized Pólya urn sequence* whose *directing measure* [in the terminology of Aldous (1983)] is a Dirichlet–Ferguson process on $(A, \mathcal{A})$ with parameter $\alpha(\cdot)/c$ [the reader is referred to Ferguson (1973), Blackwell and MacQueen (1973), Blackwell (1973) and Ferguson (1974) for definitions, proofs of the above claims and discussions of the relevance of such objects in Bayesian nonparametric statistics; see also Pitman (1996) for a rich survey of some recent developments of Pólya urn processes]. Note also that, in all cases, the law of $\mathbf{X}^{(\alpha,c)}$ is characterized by the following two facts: (i) for every $j < N$, $\mathbb{P}(X_j^{(\alpha,c)} \in dx) = \alpha(dx)/\alpha(A)$, (ii) for every $j < N$, the law of

$$\{X_{j+n}^{(\alpha,c)} : 1 \le n < N - j\}$$

under the probability measure

$$\mathbb{P}(\cdot \,|\, X_1^{(\alpha,c)} = x_1, \ldots, X_j^{(\alpha,c)} = x_j)$$

is that of a GUS of length $N - 1 - j$ and parameters $\alpha(\cdot) + \sum_{k=1,\ldots,j} \delta_{x_k}^c(\cdot)$ and $c$.

To be sure that Assumption B is satisfied and that we work with $2(N-1)$-extendible sequences, we will systematically assume that $\alpha(A) + c2(N-1) \ge 0$. For instance, in the case of extractions without replacement from a finite set of cardinality $\alpha(A) \in \mathbb{N}$, this condition is necessary and sufficient both to have $2(N-1)$-extendibility and to satisfy Assumption B. [More precisely, consider the case of extraction without replacement from a finite set $A$, and suppose that $\mathrm{Card}(A) = S > 0$, and that $S/2 < N - 1 < +\infty$. In this case, it is easy to see that every symmetric statistic of $(X_1^{(\alpha,c)}, \ldots, X_{N-1}^{(\alpha,c)})$ is contained in the space $SU_{S-N-1}(\mathbf{X}_{N-1}^{(\alpha,c)})$, i.e., the projection of any symmetric statistic on

$$\bigoplus_{k=S-N}^{N-1} SH_k(\mathbf{X}_{N-1}^{(\alpha,c)})$$

must equal zero; see Bloznelis and Götze (2001), Proposition 1, for a complete discussion of this point.]

One nice feature of GUS is that they are weakly independent, and, therefore, thanks to Theorem 6, Hoeffding-decomposable, as shown by the following:

PROPOSITION 8. *Let* $\mathbf{X}^{(\alpha,c)}$ *be a finite GUS satisfying the assumptions of this section and, for a fixed* $1 \le n < N$, *consider a symmetric* $T(\mathbf{X}_n^{(\alpha,c)}) \in$ $L^1(\mathbf{X}_n^{(\alpha,c)})$. Then, for every $m = 1, \ldots, n$ and for every $\mathbf{j}_{(n)} \in V_{N-1}(n)$ and every $\mathbf{i}_{(m)} \in V_{N-1}(m)$, the following equality holds with probability one:

$$
[T]_{n,m}^{(r)}(\mathbf{X}_{\mathbf{i}_{(m)}}) = \sum_{q=r}^{m} c^{q-r} \frac{\prod_{l=1}^{q-r}(n-r-l+1)}{\prod_{l=1}^{m-r}(\alpha(A)+c(n+l-1))} \beta_{q,m,r}(\alpha(A),c)
$$
$$
(17) \qquad \qquad \times \sum_{\mathbf{i}_{(m)} \wedge \mathbf{j}_{(n)} \subset \mathbf{j}_{(q)} \subset \mathbf{i}_{(m)}} [T]_{n,q}^{(q)}(\mathbf{X}_{\mathbf{j}_{(q)}}^{(\alpha,c)}),
$$

where $[T]_{n,0}^{(0)} = \mathbb{E}(T)$,

$$
\beta_{q,m,r}(\alpha(A),c)
$$
$$
= \begin{cases} 1, & q = m, \\ (\alpha(A)+c(m-1)) \times \cdots \times (\alpha(A)+cq), & r \leq q \leq m-1 \end{cases}
$$

and $r = r(\mathbf{i}_{(m)}, \mathbf{j}_{(n)}) = \mathrm{Card}(\mathbf{i}_{(m)} \wedge \mathbf{j}_{(n)})$, and all conventions are as before.

PROOF. To prove (17), consider a vector $\mathbf{j}_{(n)} \in V_{N-1}(n)$, as well as an index $i \notin \mathbf{j}_{(n)}$: it is easily verified that

$$
[T]_{n,1}^{(0)}(X_i^{(\alpha,c)}) = \frac{nc}{\alpha(A)+nc}[T]_{n,1}^{(1)}(X_i^{(\alpha,c)}) + \frac{\alpha(A)}{\alpha(A)+nc}[T]_{n,0}^{(0)},
$$

that gives (17) for $m = 1$. To show the general case we use once again a recurrence argument. Assume, indeed, that the result is proved for $m = 1, \ldots, k-1$: we recall that for every $\mathbf{i}_{(k)} \in V_{N-1}(k)$, for any fixed $\mathbf{x}_r = (x_1, \ldots, x_r) \in A^r$, under the probability measure

$$
\mathbb{P}[\cdot | \mathbf{X}_{\mathbf{i}_{(k)} \wedge \mathbf{j}_{(n)}}^{(\alpha,c)} = \mathbf{x}_r],
$$

where $r = r(\mathbf{i}_{(k)} \wedge \mathbf{j}_{(n)})$ is defined as in the statement, the vector $\mathbf{X}_{\mathbf{j}_{(n)} \setminus \mathbf{i}_{(k)}}^{(\alpha,c)}$ is a finite GUS of length $n - r$ and parameters $\alpha(\cdot) + \sum_{l=1,\ldots,r} \delta_{x_l}^c(\cdot)$ and $c$. Now fix $\mathbf{i}_{(k)} \in V_{N-1}(k)$ such that $r > 0$. The recurrence assumption, along with the obvious relation $(\mathbf{i}_{(k)} \setminus \mathbf{j}_{(n)}) \wedge (\mathbf{j}_{(n)} \setminus \mathbf{i}_{(k)}) = \varnothing$, implies

$$
\mathbb{E}[T(\mathbf{x}_r, \mathbf{X}_{\mathbf{j}_{(n)} \setminus \mathbf{i}_{(k)}}^{(\alpha,c)}) | \mathbf{X}_{\mathbf{i}_{(k)} \wedge \mathbf{j}_{(n)}}^{(\alpha,c)} = \mathbf{x}_r, \mathbf{X}_{\mathbf{i}_{(k)} \setminus \mathbf{j}_{(n)}}^{(\alpha,c)}]
$$
$$
= \sum_{q=0}^{k-r} c^q \frac{\prod_{l=1}^{q}(n-r-l+1)}{\prod_{l=1}^{k-r}(\alpha(A)+c(n+l-1))} \beta_{q,k-r,0}(\alpha(A)+cr,c)
$$
$$
\times \sum_{\mathbf{j}_{(q)} \subset \mathbf{i}_{(k)} \setminus \mathbf{j}_{(n)}} [T]_{n,r+q}^{(r+q)}(\mathbf{x}_r, \mathbf{X}_{\mathbf{j}_{(q)}}^{(\alpha,c)}).
$$

But

$$\beta_{q,k-r,0}(\alpha(A) + cr, c)$$
$$= \begin{cases} 1, & q = k - r, \\ (\alpha(A) + c(k-1)) \times \cdots \times (\alpha(A) + c(r+q)), & 0 \le q \le k-r-1 \end{cases}$$

and the change of variables $p = q + r$ yields immediately (17). We are left with the case $\mathbf{i}_{(k)} \wedge \mathbf{j}_{(n)} = \varnothing$: to see that the statement is still valid, fix $x_{i_k} \in A$ and $\mathbf{x}_{\mathbf{i}_{(k-1)}} = (x_{i_1}, \ldots, x_{i_{k-1}}) \in A^{k-1}$ and write, due to the recurrence assumption,

$$[T]_{n,k}^{(0)}(x_{i_k}, \mathbf{x}_{\mathbf{i}_{(k-1)}})$$
$$= \mathbb{E}[T(\mathbf{X}_{\mathbf{j}_{(n)}}^{(\alpha,c)}) | X_{i_k}^{(\alpha,c)} = x_{i_k}, \mathbf{X}_{\mathbf{i}_{(k)} \backslash i_k}^{(\alpha,c)} = \mathbf{x}_{\mathbf{i}_{(k-1)}}]$$
$$= \sum_{q=0}^{k-1} c^q \frac{\prod_{l=1}^{q}(n-l+1)}{\prod_{l=1}^{k-1}(\alpha(A) + c(n+l))} \beta_{q,k-1,0}(\alpha(A) + c, c)$$
$$\qquad \times \sum_{\mathbf{j}_{(q)} \subset \mathbf{i}_{(k)} \backslash i_k} [T]_{n,q+1}^{(q)}(\mathbf{x}_{\mathbf{j}_{(q)}}, x_{i_k})$$
$$= \sum_{q=0}^{k-1} c^q \frac{\prod_{l=1}^{q}(n-l+1)}{\prod_{l=1}^{k-1}(\alpha(A) + c(n+l))} (\alpha(A) + c(k-1))$$
$$\qquad \times \cdots \times (\alpha(A) + c(q+1))$$
$$\qquad \times \sum_{\mathbf{j}_{(q)} \subset \mathbf{i}_{(k)} \backslash i_k} \left[ \frac{c(n-q)}{\alpha(A) + cn} [T]_{n,q+1}^{(q+1)}(\mathbf{x}_{\mathbf{j}_{(q)}}, x_{i_k}) + \frac{\alpha(A) + cq}{\alpha(A) + cn} [T]_{n,q}^{(q)}(\mathbf{x}_{\mathbf{j}_{(q)}}) \right],$$

where $\mathbf{x}_{\mathbf{j}_{(q)}}$ stands for $(x_{j_1}, \ldots, x_{j_q})$, giving the desired conclusion.  $\square$

Actually, Proposition 8 yields much more than weak independence. As a matter of fact, we have the following:

COROLLARY 9. *Let $\mathbf{X}^{(\alpha,c)}$ be a GUS as in Proposition 8, and fix $1 \le n < N$ and $m < n$: if a symmetric $T$ on $A^n$ is such that $T(\mathbf{X}_n^{(\alpha,c)}) \in L^1(\mathbf{X}_n^{(\alpha,c)})$ and*

$$(18) \qquad\qquad [T]_{n,m}^{(m)}(\mathbf{X}_m^{(\alpha,c)}) = 0, \qquad \mathbb{P}\text{-}a.s.,$$

*then*

$$[T]_{n,m}^{(r)}(\mathbf{X}_m^{(\alpha,c)}) = 0, \qquad \mathbb{P}\text{-}a.s.$$

*for every $r \le m$ and such that $n + m - r < N$. In particular, if $T(\mathbf{X}_n^{(\alpha,c)}) \in L_s^2(\mathbf{X}_n^{(\alpha,c)})$ and $T$ satisfies (18), then $T(\mathbf{X}_{\mathbf{j}_{(n)}}^{(\alpha,c)}) \in U_m(\mathbf{X}_M^{(\alpha,c)})^{\perp}$ for every* $\mathbf{j}_{(n)} \in V_{N-1}(n)$ and $m \le M < N$, where $U_m(\mathbf{X}_M^{(\alpha,c)})$ denotes the direct sum of the first $m$ Hoeffding spaces associated to $\mathbf{X}_M^{(\alpha,c)}$, and the orthogonal is taken in $L^2(\mathbf{X}^{(\alpha,c)})$. This implies that $\mathbf{X}^{(\alpha,c)}$ is weakly independent.

We have also the following generalization of the calculations contained, for example, in Bloznelis and Götze [(2001), formula (2.5)].

COROLLARY 10.    *Let $\mathbf{X}^{(\alpha,c)}$ be a GUS. Take $T$ and $V$ square integrable, symmetric on $A^n$ and satisfying the hypotheses of Corollary 9 for $m = n-1$ [i.e., $T, V \in \Xi_n(\mathbf{X}^{(\alpha,c)})$]: then, for every $\mathbf{j}_{(n)}, \mathbf{i}_{(n)} \in V_{N-1}(n)$ such that $\mathrm{Card}(\mathbf{i}_{(n)} \wedge \mathbf{j}_{(n)}) = r$,*

$$(19) \qquad \begin{aligned} &\mathbb{E}[T(\mathbf{X}_{\mathbf{i}_{(n)}}^{(\alpha,c)})V(\mathbf{X}_{\mathbf{j}_{(n)}}^{(\alpha,c)})] \\ &= c^{n-r} \prod_{l=1}^{n-r} \frac{n-r-l+1}{\alpha(A)+c(n+l-1)} \mathbb{E}[T(\mathbf{X}_n^{(\alpha,c)})V(\mathbf{X}_n^{(\alpha,c)})]. \end{aligned}$$

We now want to calculate the explicit form of the Hoeffding-ANOVA decomposition for urn sequences.

5.1. *Hoeffding decompositions for GUS (statements).*    Now consider a sequence $\mathbf{X}^{(\alpha,c)}$ that is a GUS in the sense of the previous section, and fix $1 \le M < N$. Most of the subsequent results are related to the following sequence of real constants associated to the law of $\mathbf{X}^{(\alpha,c)}$:

$$(20) \quad \Phi(n,m,r,p) := c^p (m-r)_{(m-r-p)} \frac{\prod_{s=1}^{m-(r+p)}[\alpha(A)+c(r+p+s-1)]}{\prod_{s=1}^{m-r}[\alpha(A)+c(n+s-1)]},$$

where $1 \le m \le n \le M$, $0 \le r \le m$, $0 \le p \le m-r$, $\alpha(A)+c(n+m-r) > 0$, $(a)_{(b)} := a!/b!$ for $a \ge b$ and $\prod_{s=1}^0 = 1 = 0^0$ by convention, and, for $1 \le q \le m \le n \le M$,

$$(21) \qquad \Psi_M(q,n,m) := \sum_{r=0}^q \binom{q}{r} \binom{M-n}{m-r}_* \Phi(n,m,r,q-r)$$

with $\binom{a}{b}_* := \binom{a}{b} \mathbf{1}_{(a \ge b)}$. We are now in a position to state the main result of the section.

THEOREM 11.    *Under the previous notation and assumptions, fix $\alpha \in \mathcal{M}(A)$ such that $\Psi_M(q,n,q) \ne 0$ for every $n = 1, \ldots, M$ and every $1 \le q \le n$. Write also, for any $k \ge 1$,*

$$(22) \qquad \gamma_M^{(k)} := (\Psi_M(k,k,k))^{-1}$$

*with the notation introduced in* (21). *For every* $s = 1, \ldots, M-1$, *the following equality holds a.s.-$\mathbb{P}$ for any* $T \in L_s^2(\mathbf{X}_M^{(\alpha,c)})$ *with* $\mathbb{E}(T) = 0$:

$$\pi[T, SH_s](\mathbf{X}_M^{(\alpha,c)}) = \sum_{a=1}^{s} \theta_M^{(s,a)} \sum_{\mathbf{j}_{(a)} \in V_M(a)} [T]_{M,a}^{(a)}(\mathbf{X}_{\mathbf{j}_{(a)}}^{(\alpha,c)})$$

$$(23)$$

$$= \sum_{\mathbf{j}_{(s)} \in V_M(s)} \left[ \sum_{a=1}^{s} \theta_{M*}^{(s,a)} \sum_{\mathbf{j}_{(a)} \subset \mathbf{j}_{(s)}} [T]_{M,a}^{(a)}(\mathbf{X}_{\mathbf{j}_{(a)}}^{(\alpha,c)}) \right],$$

*where* $\theta_{M*}^{(k,a)} := \theta_M^{(k,a)} \binom{M-a}{k-a}^{-1}$ *and the coefficients* $\theta_M^{(k,a)}$ *are recursively defined by the set of conditions* $\{\mathbf{S}_M(k), k = 1, \ldots, M-1\}$ *given by*

$$(24) \quad \mathbf{S}_M(k) := \begin{cases} \theta_M^{(k,k)} = \gamma_M^{(k)}, \\ \sum_{i=q}^{k} \sum_{j=q}^{i} \theta_M^{(i,j)} \Psi_M(q,k,j) = 0, & q = 1, \ldots, k-1, \end{cases}$$

*and, consequently,*

$$\pi[T, SH_M](\mathbf{X}_M^{(\alpha,c)}) = \sum_{a=1}^{M} \sum_{\mathbf{j}_{(a)} \in V_M(a)} \theta_M^{(M,a)} [T]_{M,a}^{(a)}(\mathbf{X}_{\mathbf{j}_{(a)}}^{(\alpha,c)}),$$

*where* $\theta_M^{(M,a)} := -\sum_{s=a}^{M-1} \theta_M^{(s,a)}$ *for* $a = 1, \ldots, M-1$ *and* $\theta_M^{(M,M)} = \Psi_M(M,M,M)^{-1} = 1$.

Note how the above assumptions, concerning the constants $\Psi_M(\cdot, \cdot, \cdot)$, are immaterial in the case $c \geq 0$. It is also clear that Theorem 11 can be applied to noncentered symmetric statistics by considering $T' := T - \mathbb{E}(T)$.

The statement of Theorem 11 can be further refined by means of Theorem 6 and Corollaries 4 and 10. Indeed, the symmetric functionals

$$(25) \qquad \phi_T^{(s)}(\mathbf{X}_{\mathbf{j}_{(s)}}^{(\alpha,c)}) := \left[ \sum_{a=1}^{s} \theta_{M*}^{(s,a)} \sum_{\mathbf{j}_{(a)} \subset \mathbf{j}_{(s)}} [T]_{M,a}^{(a)}(\mathbf{X}_{\mathbf{j}_{(a)}}^{(\alpha,c)}) \right]$$

*defined for* $s = 1, \ldots, M$ *and for coefficients* $\theta_{M*}^{(\cdot,\cdot)}$ (*note that* $\theta_{M*}^{(M,\cdot)} = \theta_M^{(M,\cdot)}$) *as in* (23), *are uniquely determined* (*thanks to Corollary 4*), *and such that*

$$[\phi_T^{(s)}]_{s,s-1}^{(s-1)}(\mathbf{X}_{s-1}^{(\alpha,c)}) = 0, \qquad \mathbb{P}\text{-}a.s.,$$

*since the* $\mathbf{X}^{(\alpha,c)}$ *is weakly independent and, therefore, Hoeffding decomposable. Moreover, Corollary 9 yields*

$$[\phi_T^{(s)}]_{s,s-1}^{(r)}(\mathbf{X}_{s-1}^{(\alpha,c)}) = 0, \qquad \mathbb{P}\text{-}a.s.,$$

*for every* $2s - r < N + 1$.

REMARKS. (a) It is interesting that, for any fixed $M$, the coefficients $\theta_M^{(\cdot,\cdot)}$ appearing in (23) depend on the law of $\mathbf{X}^{(\alpha,c)}$ *only through the quantity* $c/\alpha(A)$, that can be interpreted as the (initial) *rate of replacement* associated to the GUS $\mathbf{X}^{(\alpha,c)}$. It follows that the Hoeffding-ANOVA decompositons of two different finite GUS with the same rate of replacement can be obtained by first calculating the $(M-1)$-ple of functions $[T]_{M,i}^{(i)}(\cdot)$, $i = 1, \ldots, M-1$, and then by implementing *exactly the same algorithm*.

(b) The above discussion shows that, not only a statistic $T \in SH_i(\mathbf{X}_M^{(\alpha,c)})$, $i \leq M$, is uniquely determined by a function $\phi_T^{(i)} \in \Xi_i(\mathbf{X})$, but also that such function can be "recovered" from $T$, through (25).

(c) Note that the recursive relation that defines the coefficients $\theta_M^{(\cdot,\cdot)}$ is different from that deduced in Zhao and Chen (1990) or Bloznelis and Götze (2001) for the case of $A$ being a finite set with cardinality $S > 2M$, endowed with the counting measure $\alpha$. However, Corollary 4 ensures that the results implied by Theorem 11 and those in the references above are equivalent. One can also compare the explicit computations of the parameters $\theta_M^{(1,1)}$, $\theta_M^{(2,1)}$ and $\theta_M^{(2,2)}$ that appear in Bloznelis and Götze [(2001), beginning of page 901] with those exhibited in Section 6.1.

Examples and applications of Theorem 11 are given in the next section and, to a much wider extent, in Peccati (2002a, b, 2003). Now we establish some relations that are used to prove Theorem 11.

5.2. *Auxiliary calculations.* Let $\mathbf{X}^{(\alpha,c)} = \mathbf{X}$ be a GUS as in the previous section (the dependence on $\alpha$ and $c$ is tacitly dropped to simplify the notation whenever there is no risk of confusion). The following result is the key step of the section:

PROPOSITION 12. *Let the previous notation prevail, and fix* $m$, $n$, $M$ *such that* $1 \leq m \leq n \leq M < N$, *as well as vectors* $\mathbf{i}_{(m)} \in V_M(m)$ *and* $\mathbf{j}_{(n)} \in V_M(n)$. *Then, for every symmetric* $T \in L^1(\mathbf{X}_M)$, *a version of* $\mathbb{E}[[T]_{M,m}^{(m)}(\mathbf{X}_{\mathbf{i}_{(m)}})|\mathbf{X}_{\mathbf{j}_{(n)}}]$ *is given by*

$$(26) \quad \sum_{p=0}^{m-r(\mathbf{i}_{(m)},\mathbf{j}_{(n)})} \sum_{\mathbf{l}_{(p)} \subset \mathbf{j}_{(n)} \setminus \mathbf{i}_{(m)}} \Phi(n,m,r,p)[T]_{M,r+p}^{(r+p)}(\mathbf{X}_{\mathbf{j}_{(n)} \wedge \mathbf{i}_{(m)}}, \mathbf{X}_{\mathbf{l}_{(p)}}),$$

*where* $r = r(\mathbf{i}_{(m)},\mathbf{j}_{(n)}) = \mathrm{Card}(\mathbf{i}_{(m)} \wedge \mathbf{j}_{(n)})$ *and the* $\Phi$'s *are given by* (20).

PROOF. By the symmetry of $T$ and of the distribution of the vector $\mathbf{X}_M$, we can assume without loss of generality that $\mathbf{j}_{(n)} = (1, \ldots, n)$, $\mathbf{j}_{(n)} \wedge \mathbf{i}_{(m)} = (1, \ldots, r(\mathbf{i}_{(m)},\mathbf{j}_{(n)}))$ and $i_{r+t} \geq n+1$ for $t = 1, \ldots, m - r(\mathbf{i}_{(m)},\mathbf{j}_{(n)})$.

Note that when $r(\mathbf{i}_{(m)}, \mathbf{j}_{(n)}) = m$, formula (26) is trivial and we shall therefore assume that $r = r(\mathbf{i}_{(m)}, \mathbf{j}_{(n)}) \in \{0, \ldots, m-1\}$. Now observe that, thanks to exchangeability, straightforward calculations yield

$$\mathbb{E}[[T]_{M,m}^{(m)}(\mathbf{X}_{\mathbf{i}_{(m)}})|\mathbf{X}_{\mathbf{j}_{(n)}}]$$

$$= \int_{A^{m-r}} \prod_{s=1}^{m-r} \frac{\alpha(dy_s) + \sum_{a=1}^{n} \delta_{X_a}^c(dy_s) + \sum_{a=1}^{s-1} \delta_{y_a}^c(dy_s)}{\alpha(A) + c(n+s-1)}$$

$$\times [T]_{M,m}^{(m)}(X_1, \ldots, X_r, y_1, \ldots, y_{m-r}) \qquad \text{a.s.-}\mathbb{P},$$

and one can, moreover, rewrite the product measure inside the integral according to the following formula:

$$\prod_{s=1}^{m-r} \left[ \alpha(dy_s) + \sum_{a=1}^{n} \delta_{X_a}^c(dy_s) + \sum_{a=1}^{s-1} \delta_{y_a}^c(dy_s) \right]$$

$$= \prod_{s=1}^{m-r} \left[ \alpha(dy_s) + \sum_{a=1}^{r} \delta_{X_a}^c(dy_s) + \sum_{a=1}^{s-1} \delta_{y_a}^c(dy_s) \right]$$

$$+ \sum_{r+1 \le l_1 \ne l_2 \ne \cdots \ne l_{m-r} \le n} \prod_{s=1}^{m-r} \delta_{X_{l_s}}^c(dy_s)$$

(27)

$$+ \sum_{p=1}^{m-r-1} \sum_{\substack{r+1 \le l_1 \ne \cdots \ne l_p \le n \\ \mathbf{h}_{(p)} \subset (1, \ldots, m-r)}} \left[ \prod_{s=1}^{p} \delta_{X_{l_s}}^c(dy_{h_s}) \right.$$

$$\times \prod_{q=1}^{m-r-p} \left[ \alpha(dy_{t_q}) + \sum_{s=1}^{p} \delta_{X_{l_s}}^c(dy_{t_q}) \right.$$

$$\left. \left. + \sum_{a=1}^{r} \delta_{X_a}^c(dy_{t_q}) + \sum_{a=1}^{q-1} \delta_{y_{t_a}}^c(dy_{t_q}) \right] \right],$$

where in the last summand we used the notation

$$\mathbf{t}_{(m-r-p)} = (t_1, \ldots, t_q, \ldots, t_{m-r-p})$$

$$:= (1, \ldots, m-r) \setminus \mathbf{h}_{(p)}.$$

Note that (27) can be easily shown for $m$, say, equal to 2, whereas the general case is proved by a standard recurrence argument. To conclude, use once again symmetry and exchangeability to have

$$\mathbb{E}[[T]_{M,m}^{(m)}(\mathbf{X}_{\mathbf{i}_{(m)}})|\mathbf{X}_{\mathbf{j}_{(n)}}]$$

$$= [T]_{M,r}^{(r)}(\mathbf{X}_{(\mathbf{j}_{(n)} \wedge \mathbf{i}_{(m)})}) \prod_{s=1}^{m-r} \frac{[\alpha(A) + c(r+s-1)]}{[\alpha(A) + c(n+s-1)]}$$

$$+ \frac{c^{m-r}(m-r)!}{\prod_{s=1}^{m-r}[\alpha(A) + c(n+s-1)]} \sum_{\mathbf{l}_{(m-r)} \subset (r+1,\dots,n)} [T]_{M,m}^{(m)}(\mathbf{X}_{\mathbf{i}_{(n)} \wedge \mathbf{j}_{(m)}}, \mathbf{X}_{\mathbf{l}_{(m-r)}})$$

$$+ \sum_{p=1}^{m-r(\mathbf{i}_{(m)},\mathbf{j}_{(n)})-1} \sum_{\mathbf{l}_{(p)} \subset \mathbf{j}_{(n)} \setminus \mathbf{i}_{(m)}} \left[ [T]_{M,r+p}^{(r+p)}(\mathbf{X}_{\mathbf{j}_{(n)} \wedge \mathbf{i}_{(m)}}, \mathbf{X}_{\mathbf{l}_{(p)}}) c^p p! \binom{m-r}{p} \right.$$

$$\left. \times \frac{\prod_{s=1}^{m-(r+p)}[\alpha(A) + c(r+p+s-1)]}{\prod_{s=1}^{m-r}[\alpha(A) + c(n+s-1)]} \right],$$

which agrees with (26) and (20). □

From Proposition 12 we obtain the following:

COROLLARY 13. *Under the assumptions of Proposition* 12, *for a fixed* $\mathbf{j}_{(n)} \in V_M(n)$, *a.s.-*$\mathbb{P}$,

$$\sum_{\mathbf{j}_{(m)} \in V_M(m)} \mathbb{E}[[T]_{M,m}^{(m)}(\mathbf{X}_{\mathbf{j}_{(m)}})|\mathbf{X}_{\mathbf{j}_{(n)}}] = \sum_{q=0}^m \sum_{\mathbf{j}_{(q)} \subset \mathbf{j}_{(n)}} \Psi_M(q,n,m) [T]_{M,q}^{(q)}(\mathbf{X}_{\mathbf{j}_{(q)}}),$$

*where the* $\Psi$'s *are defined as in* (21).

PROOF. Straightforward computation, along with Proposition 12, yields

$$\sum_{\mathbf{j}_{(m)} \in V_M(m)} \mathbb{E}[[T]_{M,m}^{(m)}(\mathbf{X}_{\mathbf{j}_{(m)}})|\mathbf{X}_{\mathbf{j}_{(n)}}]$$

$$= \sum_{q=0}^m \sum_{\mathbf{j}_{(q)} \subset \mathbf{j}_{(n)}} [T]_{M,q}^{(q)}(\mathbf{X}_{\mathbf{j}_{(q)}}) \left[ \sum_{\mathbf{j}_{(m)} \in V_M(m)} \mathbf{1}_{(\mathbf{j}_{(m)} \wedge \mathbf{j}_{(n)} \subset \mathbf{j}_{(q)})} \Phi(n,m,r,q-r) \right],$$

where $r := r(\mathbf{j}_{(m)} \wedge \mathbf{j}_{(n)}) = \text{Card}(\mathbf{j}_{(m)} \wedge \mathbf{j}_{(n)})$ as before and a simple combinatorial argument gives the desired result. □

Another consequence of Proposition 12 is:

COROLLARY 14. *Under the assumptions and notation of this paragraph, let* $T \in L_s^2(\mathbf{X}_M)$ *for some* $1 \le n \le M < N$, *and such that*

$$\mathbb{E}[T|\mathbf{X}_{\mathbf{i}_{(n-1)}}] = 0$$

*for every* $\mathbf{i}_{(n-1)} \in V_M(n-1)$. *Then,*

$$\pi[T, SH_n](\mathbf{X}_M) = \gamma_M^{(n)} \sum_{\mathbf{j}_{(n)} \in V_M(n)} [T]_{M,n}^{(n)}(\mathbf{X}_{\mathbf{j}_{(n)}})$$

with $\gamma_M^{(n)}$ defined according to (22).

Proof. We shall find a constant $q$ such that, a.s.-$\mathbb{P}$, for every $\mathbf{i}_{(n)} \in V_M(n)$,

$$[T]_{M,n}^{(n)}(\mathbf{X}_{\mathbf{i}_{(n)}}) - q \sum_{\mathbf{j}_{(n)} \in V_M(n)} \mathbb{E}[[T]_{M,n}^{(n)}(\mathbf{X}_{\mathbf{j}_{(n)}})|\mathbf{X}_{\mathbf{i}_{(n)}}] = 0.$$

But, thanks to Corollary 13 and the hypotheses in the statement, we can explicitly compute

$$\sum_{\mathbf{j}_{(n)} \in V_M(n)} \mathbb{E}[[T]_{M,n}^{(n)}(\mathbf{X}_{\mathbf{j}_{(n)}})|\mathbf{X}_{\mathbf{i}_{(n)}}] = [T]_{M,n}^{(n)}(\mathbf{X}_{\mathbf{i}_{(n)}})\Psi_M(n,n,n),$$

thus concluding the proof. $\square$

5.3. *End of the proof of Theorem* 11. To obtain the coefficients appearing in (23), just write for $s = 1, \ldots, M-1$ and a given $\mathbf{j}_{(s)} \in V_M(s)$, the a.s. condition

$$(28) \quad [T]_{M,s}^{(s)}(\mathbf{X}_{\mathbf{j}_{(s)}}) - \sum_{b=1}^{s}\sum_{a=1}^{b} \theta_M^{(b,a)} \sum_{\mathbf{j}_{(a)} \in V_M(a)} \mathbb{E}([T]_{M,a}^{(a)}(\mathbf{X}_{\mathbf{j}_{(a)}})|\mathbf{X}_{\mathbf{j}_{(s)}}) = 0$$

and observe that exchangeability and symmetry imply that if the $\theta_M$'s satisfy (28) for one $\mathbf{j}_{(s)} \in V_M(s)$, then a.s. they satisfy the same condition for every element of $V_M(s)$; but the left-hand side of (28) can be rewritten, due to Corollary 13, as

$$[T]_{M,s}^{(s)}(\mathbf{X}_{\mathbf{j}_{(s)}})(1 - \theta_M^{(s,s)}\Psi_M(s,s,s))$$

$$- \sum_{q=1}^{s-1}\sum_{\mathbf{j}_{(q)} \subset \mathbf{j}_{(s)}} [T]_{M,q}^{(q)}(\mathbf{X}_{\mathbf{j}_{(q)}}) \left[\sum_{b=q}^{s}\sum_{a=q}^{b} \theta_M^{(b,a)}\Psi_M(q,s,a)\right]$$

that implies (24). The last assertion in the statement of Theorem 11 is just plain algebra, and the proof is therefore concluded.

## 6. Examples and applications.

6.1. *Examples* (*maxima and minima*). In this section we first consider a finite GUS, noted $\mathbf{X}^{(\alpha,c)}$, for which we calculate the first two terms of the Hoeffding-ANOVA decomposition of a $T \in L^2(\mathbf{X}_M^{(\alpha,c)})$ for $M > 2$. Then, we apply such a result to the case of the simplest *order statistics* associated to a real valued finite GUS.

Now consider a finite GUS $\mathbf{X}^{(\alpha,c)}$ satisfying the assumptions of Theorem 11: it is easily seen that $\theta_M^{(1,1)} = (\alpha(A) + c)/(\alpha(A) + cM)$, $\theta_M^{(2,2)} = (\alpha(A) + 3c) \times (\alpha(A) + 2c)/(\alpha(A) + Mc)(\alpha(A) + c(M+1))$ and

$$\theta_M^{(2,1)} = -\frac{(M-1)(\alpha(A) + 3c)(\alpha(A) + c)}{(\alpha(A) + cM)(\alpha(A) + c(M+1))},$$

so that, for any symmetric and centered $T \in L^2(\mathbf{X}_M^{(\alpha,c)})$,

$$\pi[T, SH_1](\mathbf{X}_M)$$
$$= \frac{\alpha(A) + c}{\alpha(A) + cM} \sum_{i=1}^{M} [T]_{M,1}^{(1)}(X_i^{(\alpha,c)}),$$

$$\pi[T, SH_2](\mathbf{X}_M)$$
$$= \frac{(\alpha(A) + 3c)(\alpha(A) + 2c)}{(\alpha(A) + cM)(\alpha(A) + c(M+1))} \sum_{\mathbf{j}_{(2)} \in V_M(2)} [T]_{M,2}^{(2)}(\mathbf{X}_{\mathbf{j}_{(2)}}^{(\alpha,c)})$$

$$- \frac{(M-1)(\alpha(A) + 3c)(\alpha(A) + c)}{(\alpha(A) + cM)(\alpha(A) + c(M+1))} \sum_{i=1}^{M} [T]_{M,1}^{(1)}(X_i^{(\alpha,c)})$$

$$= \sum_{\mathbf{j}_{(2)} \in V_M(2)} \Bigg[ \frac{(\alpha(A) + 3c)(\alpha(A) + 2c)}{(\alpha(A) + cM)(\alpha(A) + c(M+1))} [T]_{M,2}^{(2)}(\mathbf{X}_{\mathbf{j}_{(2)}}^{(\alpha,c)})$$

$$- \frac{(\alpha(A) + 3c)(\alpha(A) + c)}{(\alpha(A) + cM)(\alpha(A) + c(M+1))} \sum_{i \subset \mathbf{j}_{(2)}} [T]_{M,1}^{(1)}(X_i^{(\alpha,c)}) \Bigg].$$

Suppose also that $A \subset \Re$: we shall compute the three quantities $\mathbb{E}(T(\mathbf{X}_M^{(\alpha,c)}))$, $[T]_{M,1}^{(1)}(z)$ and $[T]_{M,2}^{(2)}(z_1, z_2)$ associated to the symmetric statistics $T(\mathbf{X}_M^{(\alpha,c)}) = \max(X_1^{(\alpha,c)}, \dots, X_M^{(\alpha,c)})$ (the same calculations hold for the minimum), so to write the first two terms of its Hoeffding-ANOVA decomposition. In this case, it is easily seen that

$$\mathbb{E}(T(\mathbf{X}_M^{(\alpha,c)})) = \sum_{k=0}^{M-1} n(k, c, \alpha(A)) \int_{A^{M-k}} \max(x_1, \dots, x_{M-k})$$
$$(29) \qquad\qquad\qquad\qquad\qquad \times \alpha^{\otimes(M-k)}(dx_1, \dots, dx_{M-k}),$$

where

$$n(k, c, \alpha(A)) := c^k \Bigg[ \sum_{i_1=1}^{M-1} \cdots \sum_{i_k=i_{k-1}+1}^{M-1} \left( \prod_{s=1}^{k} i_s \right) \Bigg] \Bigg[ \prod_{t=1}^{M} (\alpha(A) + c(t-1))^{-1} \Bigg]$$

and, again,

$$[T]_{M,1}^{(1)}(z) = \sum_{k=0}^{M-2} n'(k, c, \alpha(A))$$

$$\times \sum_{i=0}^{M-1-k} \binom{M-k-1}{i} c^{M-k-1-i}$$

$$\times \int_{A^i} \max(x_1, \ldots, x_i, z) \alpha^{\otimes i}(dx_1, \ldots, dx_i)$$

$$= \sum_{i=0}^{M-1} \zeta(i, c, \alpha(A)) \int_{A^i} \max(x_1, \ldots, x_i, z) \alpha^{\otimes i}(dx_1, \ldots, dx_i),$$

where $n'(k, c, \alpha(A)) := n(k, c, \alpha(A)) \prod_{j=1}^{M}(\alpha(A) + c(j-1)) / \prod_{j=1}^{M-1}(\alpha(A) + cj)$
and

$$\zeta(i, c, \alpha(A)) = \sum_{k=0}^{(M-2) \wedge (M-1-i)} n'(k, c, \alpha(A)) \binom{M-k-1}{i} c^{M-k-1-i}$$

and, eventually,

$$[T]_{M,2}^{(2)}(z_1, z_2)$$

$$= \sum_{k=0}^{M-3} n''(k, c, \alpha(A))$$

$$\times \sum_{i=0}^{M-2-k} \sum_{j=0}^{i} \binom{M-k-2}{i} \binom{i}{j} c^{M-k-2-j}$$

$$\times \int_{A^j} \max(x_1, \ldots, x_j, z_1, z_2) \alpha^{\otimes j}(dx_1, \ldots, dx_j)$$

$$= \sum_{j=0}^{M-2} \zeta'(j, c, \alpha(A)) \int_{A^j} \max(x_1, \ldots, x_j, z_1, z_2) \alpha^{\otimes j}(dx_1, \ldots, dx_j)$$

with $n''(k, c, \alpha(A)) := n'(k, c, \alpha(A)) \prod_{j=1}^{M-1}(\alpha(A) + cj) / \prod_{j=1}^{M-2}(\alpha(A) + c(j+1))$, and

$$\zeta'(j, c, \alpha(A)) = \sum_{k=0}^{M-3} \sum_{i=j}^{M-2-k} n''(k, c, \alpha(A)) \binom{M-k-2}{i} \binom{i}{j} c^{M-k-2-j}$$

and, therefore, by noting $Q_{1,i}^{\alpha}(z) = \int_{A^i} \max(x_1, \ldots, x_i, z) \alpha^{\otimes i}(dx_1, \ldots, dx_i)$ and

$$Q_{2,j}^{\alpha}(z_1, z_2) = \int_{A^j} \max(x_1, \ldots, x_j, z_1, z_2) \alpha^{\otimes j}(dx_1, \ldots, dx_j),$$

we obtain

$$\pi[T, SH_1](\mathbf{X}_M^{(\alpha,c)}) = \frac{\alpha(A) + c}{\alpha(A) + cM} \sum_{i=1}^M \left( \sum_{k=0}^{M-1} \zeta(k, c, \alpha(A)) Q_{1,k}^{\alpha}(X_i^{(\alpha,c)}) - \mathbb{E}(T) \right)$$

and, finally,

$$\pi[T, SH_2](\mathbf{X}_M^{(\alpha,c)})$$

$$= \frac{(\alpha(A) + 3c)(\alpha(A) + 2c)}{(\alpha(A) + cM)(\alpha(A) + c(M+1))}$$

$$\times \sum_{\mathbf{j}_{(2)} \in V_M(2)} \left[ \left( \sum_{k=0}^{M-2} \zeta'(k, c, \alpha(A)) Q_{2,k}^{\alpha}(\mathbf{X}_{\mathbf{j}_{(2)}}^{(\alpha,c)}) - \mathbb{E}(T) \right) \right.$$

$$\left. - \frac{\alpha(A) + c}{\alpha(A) + 2c} \sum_{i \subset \mathbf{j}_{(2)}} \left( \sum_{k=0}^{M-1} \zeta(k, c, \alpha(A)) Q_{1,k}^{\alpha}(X_i^{(\alpha,c)}) - \mathbb{E}(T) \right) \right],$$

where $\mathbb{E}(T)$ is given in (29). This example shows in particular that, even if the coefficients determining the decomposition of a symmetric statistic depend exclusively on the rate of replacement associated to the GUS, the whole decomposition strongly depends on the form of the associated measure $\alpha$ [in this case through the functions $Q^{\alpha}(\cdot)$]. To conclude, observe that, for $c = 0$ and $\alpha(A) = 1$, the above calculations reduce to the usual formulae for i.i.d. random variables,

$$\pi[T, SH_1](\mathbf{X}_M^{(\alpha,0)}) = \sum_{i=1}^M [Q_{1,M-1}^{\alpha}(X_i^{(\alpha,0)}) - \mathbb{E}(T)]$$

$$\pi[T, SH_2](\mathbf{X}_M^{(\alpha,0)}) = \sum_{\mathbf{j}_{(2)} \in V_M(2)} \left[ (Q_{2,M-2}^{\alpha}(\mathbf{X}_{\mathbf{j}_{(2)}}^{(\alpha,0)}) - \mathbb{E}(T)) \right.$$

$$\left. - \sum_{i \subset \mathbf{j}_{(2)}} (Q_{1,M-1}^{\alpha}(X_i^{(\alpha,0)}) - \mathbb{E}(T)) \right].$$

6.2. *Weak copies of exchangeable sequences.* The content of this section is inspired by Föllmer, Wu and Yor (2000). Given an infinite exchangeable sequence $\mathbf{X}$ with values in a Polish space $(A, \mathcal{A})$, and $k \geq 1$, we say that a random sequence $\mathbf{Y} = \{Y_n : n \geq 1\}$ is a *k-weak copy* of $\mathbf{X}$ if, for every $\mathbf{j}_{(k)} \in V_\infty(k)$, $\mathbf{Y}_{\mathbf{j}_{(k)}} \overset{\text{law}}{=} \mathbf{X}_{\mathbf{j}_{(k)}}$. Plainly, $\mathbf{X}$ is a $k$-weak copy of itself for each $k$: however, one may wonder whether there exist $k$-weak copies of $\mathbf{X}$ for some $k$, whose law differs from that of $\mathbf{X}$. Such a problem can be solved by means of the theory developed in this paper: the next proposition shows that the answer is positive for some class of weakly independent sequences (containing infinite GUS), and that, moreover, exchangeability is preserved by any $k$-weak copy of $\mathbf{X}$ constructed by our techniques. We will note by $D(\cdot;\omega)$ the directing (probability) measure of the infinite sequence $\mathbf{X}$.

PROPOSITION 15. *Suppose that the infinite exchangeable sequence $\mathbf{X}$ is Hoeffding decomposable, and that there exists some bounded and symmetric $T$ on $A^{k+1}$ such that $T \in \Xi_{k+1}(\mathbf{X})$ and*

$$\mathbb{P}\bigg( \int_{A^{k+1}} T \, dD^{\otimes k+1} \neq 0 \bigg) > 0.$$

*Consider, moreover, the canonical space $(A^\infty, \mathcal{A}^{\otimes\infty}, \mathbb{P})$, where $\mathbb{P}$ is the law of $\mathbf{X}$. Then, there exists a random sequence $\mathbf{Y}^{(k)} = \{Y_n^{(k)} : n \geq 1\}$, with elements taking values in $(A, \mathcal{A})$ and with law $\mathbb{Q}_k(\cdot)$, that has the following properties:*

1. $\mathbb{Q}_k \ll \mathbb{P}$;
2. $\mathbb{Q}_k \neq \mathbb{P}$;
3. $\mathbf{Y}^{(k)}$ *is exchangeable;*
4. $\mathbf{Y}^{(k)}$ *is a $k$-weak copy of $\mathbf{X}$.*

*Moreover, for every $\eta > 0$, there exists $\mathbb{Q}_{k,\eta}$ satisfying points 1 to 4 above and such that*

$$\bigg\| 1 - \frac{d\mathbb{Q}_{k,\eta}}{d\mathbb{P}} \bigg\|_\infty < \eta.$$

PROOF. Call $\widetilde{\mathbf{X}} = \{\widetilde{X}_n : n \geq 1\}$ the canonical projection of the space $(A^\infty, \mathcal{A}^{\otimes\infty})$ to itself so that, if we endow $(A^\infty, \mathcal{A}^{\otimes\infty})$ with a probability measure $\mathbb{P}$, $\widetilde{\mathbf{X}}$ becomes a random element with law $\mathbb{P}$. Now, it is immediate that under the probability measure $\mathbb{Q}_k$ given by

$$d\mathbb{Q}_k = \bigg( 1 + \int_{A^{k+1}} T \, dD^{\otimes k+1} \bigg) d\mathbb{P},$$

$\widetilde{\mathbf{X}}$ satisfies points 1 to 3 in the statement: moreover, for every $\mathbf{j}_{(k)} = (j_1, \ldots, j_k) \in V_\infty(k)$,

$$\mathbb{Q}_k(\widetilde{X}_{j_1} \in B_1, \ldots, \widetilde{X}_{j_k} \in B_k) = \mathbb{E}\bigg[ \bigg( 1 + \int_{A^{k+1}} T \, dD^{\otimes k+1} \bigg) \mathbb{1}_{(\widetilde{X}_{j_1} \in B_1, \ldots, \widetilde{X}_{j_k} \in B_k)} \bigg]$$

$$= \mathbb{P}(\widetilde{X}_{j_1} \in B_1, \ldots, \widetilde{X}_{j_k} \in B_k)$$

due to the weak independence of $\mathbf{X}$ as well as the following relation, that is a consequence of de Finetti's theorem and of the fact that $T$ is bounded,

$$\lim_{M \to +\infty} \binom{M}{k+1}^{-1} \sum_{\mathbf{j}_{(k+1)} \in V_M(k+1)} T(\mathbf{X}_{\mathbf{j}_{(k+1)}}) = \int_{A^{k+1}} T \, dD^{\otimes k+1}, \qquad \mathbb{P}\text{-a.s.,}$$

yielding point 4. The last assertion is an easy consequence of the above discussion. $\qquad \square$

To eventually construct such a $T$ for a GUS $\mathbf{X}$ of parameters $\alpha$ and $c \geq 0$, maintain the notation of the proof of the above proposition, set $\mathbb{P} = \mathbb{P}^{(\alpha,c)}$, that is, the law of $\mathbf{X}$, and take a bounded and symmetric statistic $V(\widetilde{\mathbf{X}}_{k+1})$. Theorem 11 implies that one can choose $V$ such that the functional $\pi[V, SH_{k+1}](\widetilde{\mathbf{X}}_{k+1})$ is not only symmetric and different from zero, but also a.s.-$\mathbb{P}^{(\alpha,c)}$ equal to a finite linear combination of conditional expectations of $V$. It follows that for any $\eta \in (0,1)$, there exists $\varepsilon > 0$ such that $\varepsilon|\pi[V, SH_{k+1}]| < \eta$, $\mathbb{P}^{(\alpha,c)}$-a.s. It is shown in Peccati (2002a), that in this case

$$\mathbb{P}\left(\int_{A^{k+1}} \pi[V, SH_{k+1}] \, dD^{\otimes k+1} \neq 0\right) > 0$$

so that it is sufficient to take $T = \varepsilon\pi[V, H_{k+1}]$.

6.3. *Covariance analysis.* A standard combinatorial argument yields the following result that shows how the covariance of two centered and symmetric statistics can be decomposed by means of the functions $\phi^{(i)}$ defined in (25).

PROPOSITION 16 (Covariance decomposition). *Under the assumptions of Theorem* 11, *let $T$ and $Z$ be two centered elements of $L_s^2(\mathbf{X}_M^{(\alpha,c)})$, $1 \leq M < N$, and let the functions $\phi_T^{(s)}$ and $\phi_Z^{(s)}$, $s = 1, \ldots, M$, be defined by* (25). *Then*

$$\mathbb{E}[TZ] = \sum_{s=1}^{M} J_M(s, c, \alpha(A))\mathbb{E}[\phi_T^{(s)}(\mathbf{X}_s)\phi_Z^{(s)}(\mathbf{X}_s)],$$

*where*

$$J_M(s, c, \alpha(A)) := \binom{M}{s}\sum_{p=0}^{s}\binom{s}{p}\binom{M-s}{s-p}_* c^{s-p}\prod_{l=1}^{s-p}\frac{s-p-l+1}{\alpha(A)+c(s+l-1)}.$$

**Acknowledgments.** I am grateful to Professor D. M. Cifarelli and Professor P. Muliere for introducing me to generalized Pólya sequences. Section 4 has benefited from a conversation with Professor J. Pitman, in the occasion of the Saint Flour summer school of July 2002. The expression "weakly independent random variables" emerged during an inspiring discussion with Professor C. Houdré, to whom I wish to express my gratitude.

This paper is dedicated to the memory of Marie Bauer.

Laboratoire de Statistique Théorique et Appliquée
Université de Paris VI
4 Place Jussieu
75252 Paris Cedex 05
France

Istituto di Metodi Quantitativi dell'Università L. Bocconi 25
via Sarfatti
20136 Milan
Italy
e-mail: giovanni.peccati@libero.it